\documentclass[pdflatex]{article}
\usepackage{latexsym}
\usepackage{amssymb}
\usepackage{amsmath}
\usepackage{amsfonts}
\usepackage{mathrsfs}

\usepackage{color}

\newcommand{\T}{\mathcal{T}}

\newcommand{\R}{\mathbb{R}}

\newcommand {\ep}  {\varepsilon}
\newcommand {\proof} {\noindent {\bf Proof}. }

\newcommand {\p}   {\partial}
\newcommand {\f}   {\frac}
\newcommand{\beq}{\begin{equation}}
\newcommand{\eeq}{\end{equation}}
\newcommand{\qed}{{ \hfill
                       {\unskip\kern 6pt\penalty 500
                       \raise -2pt\hbox{\vrule\vbox to 6pt{\hrule width 6pt
                       \vfill\hrule}\vrule} \par}   }}
\usepackage{latexsym}
\usepackage{epsfig, graphicx}

\newcommand{\fer}[1]{(\ref{#1})}

\newcommand{\commentout}[1]{}

\newcommand {\ud}{\, \mathrm d}
\newcommand{\dis}{\displaystyle}

\newcommand{\bea} {\begin{array}{rl}}
\newcommand{\eea} {\end{array}}
\newcommand{\bepa}{\left\{ \begin{array}{l}}
\newcommand{\eepa} {\end{array}\right.}
\newtheorem{theorem}{Theorem}%[section]
\newtheorem{lemma}{Lemma}%[section]
\newtheorem{definition}{Definition}%[section]
\newtheorem{corollary}{Corollary}%[section]
                       \title{Scaling limit of a discrete prion dynamics model}
\author{Marie Doumic  \thanks{INRIA Research Centre Paris--Rocquencourt, Project-Team BANG, Domaine de Voluceau, BP 105, 781153 Rocquencourt, France; emails: marie.doumic@inria.fr, thomas.lepoutre@inria.fr} \and Thierry Goudon \thanks{INRIA Research Centre Lille-Nord Europe, 
Project-Team SIMPAF,
Park Plazza, 40 avenue Halley, F-Villeneuve d'Ascq cedex, France ; email: thierry.goudon@inria.fr} \and Thomas Lepoutre\thanks{UPMC Univ Paris 06, UMR 7598, Laboratoire Jacques-Louis Lions, F-75005, Paris, France}
\thanks{CNRS, UMR 7598, Laboratoire Jacques-Louis Lions, F-75005, Paris, France} \footnotemark[1]}
\begin{document}
\maketitle

\noindent{\bf Abstract}
This paper investigates the connection between discrete and continuous models describing
prion proliferation. The scaling
parameters are interpreted on 
biological grounds and we establish rigorous convergence statements. We also discuss, based on the asymptotic analysis, relevant boundary conditions that can be used to complete the continuous model.
 
 \


\

\noindent{\bf AMS Class. No.} 35B45, 45K05, 92D25

\

\section*{Introduction}

The modelling of intracellular prion infection has been dramatically improved in the past few years according to recent progress in molecular biology of this pathology. Relevant models have been designed to investigate the conversion of the normal monomeric form of the protein (denoted PrPc) into the infectious polymeric form (denoted PrPsc) according to the auto-catalytic process: 
$$\mathrm{PrPc}+\mathrm{PrPsc} \longrightarrow 2\mathrm{PrPsc},$$ 
in fibrillar aggregation of the protein. These models are based on linear growth of PrPsc polymers via an autocatalytic process \cite{Eigen1996}. 

The seminal paper by Masel et al. \cite{Masel} proposed a discrete model where the prion population is described by its distribution with respect to the size of polymer aggregates. The model is an infinite-dimensional system of Ordinary Differential Equations, taking into account nucleated transconformation and polymerization, fragmentation and degradation of the polymers, as well as production of PrPc by the cells. This model consists in an aggregation fragmentation discrete model. In full generality, it writes as follows:
%%%%%%%%%%%%%%%%%%%%%%%%%%%%%%
%%%%%%%%%%%%%%%%%%%%%%%%%%%%%%EQUATION DISCRETE
%%%%%%%%%%%%%%%%%%%%%%%%%%%%%%
\begin{equation}\left\lbrace\begin{array}{l}\label{eq:discret}
\dfrac{\ud v}{\ud t}=\lambda-\gamma v -v\dis\sum_{i=n_0}^\infty \tau_i u_i +2\dis\sum_{j\geq n_0}\dis\sum_{i<n_0}i k_{i,j}\beta_j u_j,\\[0.3cm]
\dfrac{\ud u_i}{\ud t}=-\mu_i u_i-\beta_iu_i-v(\tau_iu_i-\tau_{i-1}u_{i-1})+2\dis\sum_{j>i}\beta_jk_{i,j}u_j,\quad \text{for } i\geq n_0 
\end{array}\right.
\end{equation}
(with the convention $\tau_{n_0-1}u_{n_0-1}=0$).
%%%%%%%%%%%%%%%%%%%%%%%%%%%%%%
%%%%%%%%%%%%%%%%%%%%%%%%%%%%%%EQUATION DISCRETE
%%%%%%%%%%%%%%%%%%%%%%%%%%%%%%
Here $v$ represents the quantity of healthy monomers (PrPc), $u_i$ the quantity of infectious polymers (PrPsc) of size $i$, \emph{i.e.} formed by the fibrillar aggregation of $i$ monomers. We thus have $i\geqslant n_0 \geqslant 2,$ where $n_0$ represents the minimal size for polymers: smaller polymers are considered to be unstable and are immediately degraded into monomers, as the last term of Equation (\ref{eq:discret}) for $v$ expresses. $\gamma$ and $\mu_i$ are the degradation rates respectively of monomers and polymers of size $i.$ $\lambda$ is a source term: the basal synthesis rate of PrPc. $\beta_i$ is the fragmentation rate of a polymer of size $i,$ and $k_{j,i}$ is the repartition function for a polymer of size $i$ dividing into two polymers of smaller sizes $j$ and $i-j.$ Finally, $v\tau_i$ is the aggregation speed of polymers of size $i,$  which is supposed to depend both on the available quantity of monomers $v$ and on a specific aggregation ability $\tau_i$ of polymers of size $i.$

In the original model \cite{Masel}, the degradation rate of polymers $\mu_i$ and the  aggregation rate  $\tau_i$ were assumed to be  independent of the size $i,$ the fragmentation rate satisfied $\beta_j=\beta (j-1)$ for a constant $\beta$ and $k_{i,j}$ was a uniform repartition over $\{1, \dots j-1\}$, i.e., $k_{i,j}=\frac{1}{j-1}.$ These laws express that all polymers behave in the same way, and that any joint point of any polymer has the same probability to break. It allowed the authors to close the system into an ODE system of three equations, which is quite simple to analyze. However, following recent experimental results such as in \cite{Silveira}, and their mathematical analysis in \cite{Calvez, Calvez2}, we prefer here to consider variable coefficients in their full generality. 
Following the ideas of \cite{DaCosta}, we can consider, under reasonable growth assumptions on the coefficients, the so called \emph{admissible} solutions, i.e., solutions obtained by taking the limit of truncated systems (see Appendix \ref{app:admis}). 

Recent work by Greer \emph{et al.} analyzed this process in a continuous setting \cite{Greer}. They proposed a Partial Differential Equation to render out the above-mentioned polymerization/fragmentation process. It writes
%%%%%%%%%%%%%%%%%%%%%%%%%%%%%%
%%%%%%%%%%%%%%%%%%%%%%%%%%%%%%EQUATION CONTINUE
%%%%%%%%%%%%%%%%%%%%%%%%%%%%%%
\begin{eqnarray}
\frac{\ud V}{\ud t}=\lambda-\gamma V -V\int_{x_0}^\infty \tau(x)U(t,x)\ud x +2\int_{x=x_0}^\infty\int_{y=0}^{x_0}yk(y,x)\beta(x)U(t,x)\ud x\ud y,\label{eq:ODE}\\[0.3cm]
\frac{\partial  U}{\partial t}=-\mu(x) U(t,x)-\beta(x)U(t,x)-V\f{\p}{\p x}(\tau U)+2\int_x^\infty\beta(y)k(x,y)U(t,y)\ud y.\label{eq:PDEstrong}
\end{eqnarray}
%%%%%%%%%%%%%%%%%%%%%%%%%%%%%%
%%%%%%%%%%%%%%%%%%%%%%%%%%%%%%EQUATION CONTINUE
%%%%%%%%%%%%%%%%%%%%%%%%%%%%%%
The coefficients of the continuous model \eqref{eq:ODE}\eqref{eq:PDEstrong} have the same meaning than those of the discrete one (\ref{eq:discret}); however, some questions about their scaling remain, and in particular about the exact biological interpretation of the variable $x.$ 

\

The aim of this article is to investigate the link between System \eqref{eq:discret} and System \eqref{eq:ODE}\eqref{eq:PDEstrong}. We discuss in details 
the convenient mathematical assumptions under which we can ensure that the continuous system is the limit of the discrete one and
we establish 
rigorously the convergence statement. We also want to discuss possible biological interpretations of our asymptotic analysis, and see how our work can help to define a proper boundary condition at $x=x_0$ for System \eqref{eq:ODE}\eqref{eq:PDEstrong}.
Indeed, Eq. \eqref{eq:PDEstrong} holds in the domain $x> x_0$ and, due to the convection term, at least when $V(t)\tau(x_0)>0$ 
a boundary condition is necessary to complete the problem.

\

In Section \ref{sec:basic}, we first recall general properties and previous results on the considered equations. In Section \ref{sec:notations}, we rescale the equations in order to make a small parameter $\ep$ appear, and state the main result: the asymptotic convergence of the rescaled discrete system towards the continuous equations. Section \ref{sec:moments} is devoted to its proof, based on moments \emph{a priori} estimates. Sections \ref{sec:bound} and \ref{sec:discuss} discuss how  these results can be interpreted on physical grounds. We also comment the issue of the  boundary condition for the continuous model.

%%%%%%%%%%%%%%%%%%%%%%%%%%%%%%%%%%%%%%%%%%%%%%%%%%%%%%%%%%%%%%%%%
%%%%%%%%%%%%%%%%%%%%%%%%%%%%%%%%%%%%%%%%%%%%%%%%%%%%%%%%%%%%%%%%
\section{Basic properties of the equations}\label{sec:basic}
%%%%%%%%%%%%%%%%%%%%%%%%%%%%%%%%%%%%%%%%%%%%%%%%%%%%%%%%%%%%%%%%
%%%%%%%%%%%%%%%%%%%%%%%%%%%%%%%%%%%%%%%%%%%%%%%%%%%%%%%%%%%%%%%

All the considered coefficients are nonnegative. We need some structural hypothesis on $k$ and $k_{j,i}$ to make sense. 
Obviously, the hypothesis take into account that a polymer can only break into smaller pieces.
We also impose symmetry since a given polymer of size $y$ breaks equally into two polymers of size respectively $x$ and $y-x$.
Summarizing, we have
\begin{equation}\label{as:k_struc}\begin{array}{ll}%\tag{K1}
 k_{i,j}\geq 0,\qquad & k(x,y)\geq 0,
\\
k_{j,i}=0 \textrm{ for $j\geq i$}\qquad &
k(x,y)=0\textrm{ for $x> y$},
\end{array}
\end{equation}
\begin{equation}\label{as:k_symmetric}%\tag{K1}
k_{j,i}=k_{i-j,i},\qquad k(x,y)=k(y-x,y),
\end{equation}
\begin{equation}\label{as:k_proba}%\tag{K2}
\sum_{i=1}^{j-1} k_{i,j}=1,\qquad \int_0^y k(x,y)\ud x=1.
\end{equation}
(Note that \eqref{as:k_struc} and \eqref{as:k_proba} imply that $0\leq k_{i,j}\leq 1$.)
Classically, these two conditions lead to a third one, expressing mass conservation through the fragmentation process:
\begin{equation}\label{as:k_mass}%\tag{K3}
2\sum_{i=1}^{j-1} i k_{i,j}=j,\qquad 2\int_0^y xk(x,y)\ud x=y.
\end{equation}

The discrete equation belongs to the family of coagulation-fragmentation models (see \cite{Ball-Carr-DC},\cite{Ball-Becker-Doring}). 
Adapting the work of \cite{Ball-Carr-DC,Ball-Becker-Doring} to this system, we obtain the following result. It is not optimal but sufficient for our study.

\begin{theorem}\label{th:exist:discret}

Let $k_{i.j}$ satisfy Assumptions   \fer{as:k_struc}--\fer{as:k_proba}.
We assume the following growth assumptions on the coefficients
\begin{equation}\label{growthcond}\left\{\begin{array}{l}
\textrm{There exists $K>0$, $\alpha\geq0$, $m\geq 0$ and $0\leq\theta\leq 1$ such that }
\\
0\leq\beta_i\leq Ki^\alpha
\qquad
0\leq\mu_i\leq Ki^m,\qquad
0\leq\tau_i\leq Ki^\theta.
\end{array}\right.\end{equation}
The initial data $v^0\geq 0,u_i^0\geq 0$ satisfies, for $\sigma=\max(1+m,1+\theta,\alpha)$ 
$$%v^0+\sum_{i=n_0}^\infty i u_i^0<+\infty,\qquad 
\sum_{i=n_0}^\infty i^\sigma u_i^0<+\infty.$$ 
Then there exists a unique global solution to \fer{eq:discret} which satisfies
for all $t\geq 0$
\beq v(t)+\sum_{i=n_0}^\infty iu_i(t)=v^0+\sum_{i=n_0}^\infty iu_i^0+\lambda t-\int_0^t \gamma v(s)\ud s-\int_0^t\sum_{i=n_0}^\infty i\mu_iu_i(s)\ud s.\label{eq:mass_discret}\eeq
\end{theorem}
A sketch of the proof is given in Appendix \ref{app:admis}.
Let us introduce the quantity 
\beq\label{def:rho}
\rho(t)=v(t)+\sum_{n_0}^\infty iu_i(t),
\eeq
which is the total number of monomers in the population. Equation
\fer{eq:mass_discret} is a mass balance equation, which can be written as 
\beq \label{eq:balance_dis}
\frac{\ud}{\ud t}\rho=\lambda-\gamma v(t)-\sum_{i=n_0}^\infty i\mu_i u_i(t).
\eeq 
Similarly for the continuous model we define 
$$\varrho(t)=V(t)+\int_{x_0}^\infty xU(t,x)\ud x.$$
The analogue of
\eqref{eq:balance_dis} would be
\beq \label{eq:balance}
\varrho(t)-\varrho(0)=\lambda t-\dis\int_0^t\gamma V(s)\ud s-\int_0^t\int_{x_0}^\infty x\mu(x)U(t,x)\, \ud x.
\eeq 
In fact, the argument to 
deduce \eqref{eq:balance}
from the system \eqref{eq:ODE}\eqref{eq:PDEstrong} is two--fold: it relies both 
 on the boundary condition on $\{x=x_0\}$ for \eqref{eq:PDEstrong}
and on the integrability properties of the fragmentation term
\[
x\times \Big(2\displaystyle\int_x^\infty \beta(y)k(x,y)U(t,y)\ud y
-\beta(x) U(t,x)
\Big),
 \]
the integral of which has to be combined to \eqref{eq:ODE}
by virtue of \eqref{as:k_mass}.
The question is actually quite deep, as it is already revealed by the case where $\mu=0$, $\tau=0$ and $x_0=0$.
In this situation it can be shown that \eqref{eq:PDEstrong} admits solutions that do not satisfy the conservation law: $\int_0^\infty xU(t,x)\ud x=\int_0^\infty xU(0,x)\ud x$, see \cite{DubSte}.
Hence,  \eqref{eq:balance} has to be incorporated 
in the model as a constraint to select the physically relevant solution, as suggested in 
\cite{DubSte} and \cite{Well-posed}.
Nevertheless, the integrability of 
the fragmentation term is not a big deal since it can be obtained by imposing 
boundedness of a large enough moment of the initial data
as it will be clear in the discussion below and as it appeared in \cite{DubSte, Well-posed}.
More interesting is how to interpret this in terms of boundary conditions; we shall  discuss 
the point in Section \ref{sec:bound}. 
(Note that in \cite {Well-posed} the problem is completed with the boundary condition $U(t,x_0)=0$
while $x_0>0$, $\tau(x_0)>0$.)
According to \cite{DubSte, Well-posed}
we adopt the following definition.

\begin{definition}\label{def_sol}
We say that the pair $(U,V)$ is a ``monomer preserving weak solution of the prion proliferation equations'' 
with initial data $(U_0,V_0)$
if it satisfies \eqref{eq:ODE} and if
for any $\varphi\in {{\cal C}}_c^\infty( (x_0,\infty))$,
we have 
\beq\label{eq:PDEweak}
\begin{array}{l}
\dis\int_0^\infty U(t,x)\varphi(x)\ud x-\dis\int_0^\infty U_0(x)\varphi(x)\ud x
\\
[0.3cm]
=-
\dis\int_0^t\int_0^\infty \mu(x)U(s,x)\varphi(x)\ud x\ud s-
\dis\int_0^t\int_0^\infty \beta(x)U(s,x)\varphi(x)\ud x\ud s\\
[0.3cm]
+\dis\int_0^t V(s)\int_0^\infty \tau(x)U(s,x)\partial_x\varphi(x)\ud x\ud s
+
2\dis\int_0^t\dis\int_{x_0}^\infty \beta(y)U(s,y)\int_{x_0}^y k(x,y)\varphi(x)\ud x\ud y\ud s,\\[0.3cm]
\end{array}
\eeq
and 
\beq\label{eq:mass_continu}\begin{array}{lll}
V(t)+\dis\int_{x_0}^\infty xU(t,x)\ud x
&=&
V_0+\dis\int_{x_0}^\infty x U_0(x)\ud x
\\&&
+\lambda t-
\dis\int_0^t\gamma V(s)\ud s-\dis\int_0^t\int_{x_0}^\infty x \mu(x)U(s,x)\ud x\ud s.
\end{array}\eeq
\end{definition}

A break is necessary to discuss the functional framework to be used in Definition \ref{def_sol}.
We start with a set up of a few notation.
We denote by ${\cal M}^1(X)$ the set of bounded Radon measures 
on a borelian set  $X\subset \mathbb R$;  ${\cal M}^1_+(X)$ stands for the positive cone in 
${\cal M}^1(X)$. 
The space ${\cal M}^1(X)$ identifies as the dual of the space ${\cal{C}}_0(X)$ of continuous functions vanishing at infinity in $X$,\footnote{$\phi\in  {\cal C}_0(X)$ means hat $\phi$ is continuous and for any $\eta>0$, there exists a compact set $K\subset X$ such that $\sup_{X\setminus K}|\phi(x)|\leq \eta$.
We denote  ${\cal C}_c(X)$  the space of continuous functions with compact support in $X$.}
 endowed with the supremum norm,   see \cite{Malliavin}.
Given an interval $I\subset \mathbb R$, we consider measure valued functions
$W:y\in I\mapsto W(y)\in  {\cal M}^1(X)$. Denoting $W(y,x)=W(y)(x),$ we say that 
$W\in  {\cal C}(I;{\cal M}^1(X)-\mathrm{weak}-\star)$,
if, for any $\varphi \in {\cal C}_0(X)$, the function
$y\mapsto \int_X \varphi(x)W(y,x)\ \ud x$ is continuous on $I$.
We are thus led to assume 
\[U\in {\cal C}([0,T];{\cal M}^1_+([0,\infty) )-\mathrm{weak}-\star),\qquad
V\in {\cal C}([0,T]),\]
with furthermore 
\[\mathrm{supp}\big(U(t,.)\big)\subset[x_0,\infty),\qquad \dis\int_{x_0}^\infty xU(t,x)\ \ud x<\infty,\]
which corresponds to the physical meaning of the unknowns.
Hence, formula \eqref{eq:PDEweak} makes sense for continuous coefficients
\[\mu,\quad\beta,\quad\tau\in {\cal C}([x_0,\infty)).\]
Concerning the fragmentation kernel, 
it suffices to suppose 
\[y\mapsto k(\cdot,y)\in {\cal C}([x_0,\infty);{\cal M}^1_+([0,\infty) )-\mathrm{weak}-\star).\]
%As said above it should be discussed whether \eqref{eq:balance}
%is connected to the evolution equation \eqref{eq:ODE}. We shall go back to this question in Section %\ref{sec:bound}.

\section{Main result}\label{sec:notations}

\subsection{Notations and rescaled equations}

We first rewrite system (\ref{eq:discret}) in a dimensionless form, as done for instance in \cite{Thierry} (see also \cite{MischlerLaurencot}). We summarize here all the absolute constants that we will need in the sequel:
\begin{itemize}
\item $T$  characteristic time,
\item $\cal U$ characteristic value for the concentration  of polymers $u_i$,
\item ${\cal V}$  characteristic value for the concentration of monomers $v$,
\item $\T$  characteristic value for the polymerisation rate $\tau_i$,
\item $B$ characteristic value for the fragmentation frequency $\beta_i$,
\item $d_0$  characteristic value for the degradation frequency of polymers $\mu_i$,
\item $\Gamma$ characteristic value for the degradation frequency of monomers $\gamma$,
%\item $M$  characteristic value for the total mass,
\item $L$  characteristic value for the source term $\lambda,$
%\item $m_1$ mass of one monomer.
\end{itemize}
The dimensionless quantities are defined by
$$
\begin{array}{lllll}\bar t =\dis\frac{t}{T},&
\bar v (\bar t ) = \dfrac{v(\bar t T)}{\cal V},&
\bar u_i (\bar t ) = \dfrac{u_i(\bar t T)}{\cal U},&
\bar \beta_i = \dfrac{\beta_i}{B},&
\bar \tau_i = \dfrac{\tau_i}{\T},\\[0.3cm]
\bar \mu_i =  \dfrac{\mu_i}{d_0},&
\bar \lambda =\dfrac{\lambda}{L},&%\\
\bar \gamma = \dfrac{\gamma}{\Gamma}.&
%\bar \rho =\dfrac{\rho}{M}.
\end{array}
$$
We remind that  $k_{i,j}$ is already dimensionless. 
The following dimensionless parameters appear
\begin{equation}
\left\lbrace\begin{array}{llll}
a=\dis\f{LT}{\cal V},&  b=BT,& c=\Gamma T, &  d=d_0T
\\
s=\dis\f{\cal U}{\cal V},& \nu =T\T {\cal V}. %& \eta=\frac{M}{m_1 {\cal V}}
&
\end{array}\right.
\end{equation}
Omitting the overlines, the equation becomes 
\begin{equation}\left\lbrace\begin{array}{l}\label{eq:discret_dimensionless}
\dfrac{\ud  v}{\ud t}=a \lambda-c \gamma v -\nu s v\sum\tau_i u_i +2bs\dis\sum_{j\geq n_0}\dis\sum_{i<n_0}i k_{i,j}\beta_j u_j,\\[0.3cm]
\dfrac{\ud  u_i}{\ud t}=-d \mu_i u_i-b\beta_iu_i-\nu v(\tau_iu_i-\tau_{i-1}u_{i-1})+2b\sum_{j>i}\beta_jk_{i,j}u_j,\qquad \text{for } i\geq n_0.
\end{array}\right.
\end{equation}
%We notice first that the parameter $T$ is a multiplier for each dimensionless coefficient appearing on the equations; hence we can without loss of generality suppose
%$$T=1.$$

The definition \eqref{def:rho} of the total mass  in dimensionless form becomes
\begin{equation}\label{dimlessmass}
v+s\dis\sum_{i=n_0}^\infty iu_i =%\eta 
\rho.
\end{equation}

The rationale motivating the scaling can be explained as follows.
Let $0<\ep\ll 1$ be a parameter  intended to tend to 0.
We pass from the discrete model  to the continuous model 
by 
associating to the $u_i$'s a stepwise constant function, constant  on each interval
$(\ep i,\ep(i+1))$.
Then sums over 
the index $i$ will be interpreted as Riemann sums which are expected to tend to integrals in the continuum limit while finite differences will give rise to derivatives.
Having in mind the case of homogeneous division and polymerization  rates $\beta(x)=x^\alpha$, $\tau(x)=x^\theta$, $\mu(x)=x^m$, which 
generalizes the constant-coefficient case proposed by \cite{Greer},  and their discrete analogue $\beta_i=i^\alpha$, $\tau_i=i^{\theta}$, $\mu_i=i^m$,
we shall assume that the rescaled coefficients $\beta_i,\mu_i,\tau_i$ fulfill  \eqref{growthcond}.
Therefore, we are led to set
$$%\eta =1, 
\qquad s=\varepsilon^2, $$
so that \eqref{dimlessmass} becomes
\begin{equation}\label{dimlessmasseps}
v+\ep\dis\sum_{i=n_0}^\infty \ep i\ u_i = \rho,
\end{equation}
to be compared to the definition of $\varrho$ in \eqref{eq:balance}.
This scaling means that the typical concentration 
of any aggregate 
with size $i>n_0$ is small compared to the monomers concentration, but the total mass of the aggregates is in the order of the mass of monomers.
Next, we set
\[
a=1, \quad b=\ep^\alpha,\quad
c=1,\quad d=\ep^m,\quad \nu=\ep^{\theta-1}.\]
The  rescaled equations read
%%%%%%%%%%%%%%%%%%%%%%%%%%%%%%%%%%%%%%%%%%%%%%%%%%%%%DISCRET AVEC EPSILON%%%%%%%%%%%%%%%%%%%%%%%%%%%%%%%%%%%%%%%%%%%%%%%%%%
%
%
\begin{equation}\left\lbrace\begin{array}{l}\label{eq:discrete_epsilon}
\dfrac{\ud v}{\ud t}=\lambda-\gamma v -\ep^{\theta +1} v\sum\tau_i u_i +2\varepsilon^{2+\alpha}\dis\sum_{i\geq n_0}\sum_{j<n_0}jk_{j,i}\beta_iu_i,\\[0.3cm]
\dfrac{\ud u_i}{\ud t}=-\ep^m\mu_i u_i-\ep^\alpha\beta_iu_i-\ep^{\theta -1} v(\tau_iu_i-\tau_{i-1}u_{i-1})+2\ep^\alpha\dis \sum_{j>i}\beta_jk_{i,j}u_j,\qquad \text{for } i\geq n_0.
\end{array}\right.
\end{equation}
Eventually, the threshold value $n_0$ also depends on the scaling parameter 
and we assume
\begin{equation}\label{as:n0}
 \dis\lim_{\ep\to 0} \ep n_0(\ep) =x_0 \geq 0. 
\end{equation}
This choice is discussed in Section \ref{sec:n0}.
\\

Equation \eqref{eq:discrete_epsilon} is completed by  an initial data 
$(u^{0,\ep}_i,v^{0,\ep})$ verifying, for some constants $M_0,\rho^0, M_{1+\sigma}$ independent of $\ep :$
\begin{equation}\label{cieps}
\left\{\begin{array}{l}
v^{0,\ep}+\ep^2\dis\sum_{i=n_0(\ep)}^\infty i u_i^{0,\ep}=\rho^0<+\infty,
\\
\ep\dis\sum_{i=n_0(\ep)}^\infty u_i^{0,\ep}\leq M_0<+\infty,\\
\ep^{2+\sigma}\dis\sum_{i=n_0(\ep)}^\infty i^{1+\sigma}u_i^{0,\ep}\leq M_{1+\sigma}<+\infty,
\qquad 1+\sigma>\max(1,\alpha,1+m,1+\theta).\end{array}\right.\end{equation}
For any $0<T<\infty$, Theorem \ref{th:exist:discret}
guarantees the existence of a solution $(u^\ep_i,v^\ep)$ of \eqref{eq:discrete_epsilon}.
Let us set
 \[\chi_i^\ep(x)=\chi_{[i\ep,(i+1)\ep)}(x)\]
  with $\chi_A$ the indicator function of a set $A$. 
  We introduce the piecewise constant function   
\[u^\ep(t,x):=\dis\sum_{i=n_0(\ep)}^\infty u^\ep_i(t)\chi_i^\ep(x)\]
On the same token, we associate  the following functions to the coefficients
\[\begin{array}{l}
 k^\ep(x,y):=\dis\sum_{i=0}^\infty\sum_{j=0}^\infty \f{k_{i,j}}{\ep}\chi_i^\ep(x)\chi_j^\ep(y),\\
\mu^\ep(x):=\dis\sum_{i=n_0(\ep)}^\infty\ep^m \mu_i\chi_i^\ep(x),
\qquad
\beta^\ep(x):=\dis\sum_{i=n_0(\ep)}^\infty \ep^{\alpha}\beta_i\chi_i^\ep(x),
\qquad
\tau^\ep(x):=\dis\sum_{i=n_0(\ep)}^\infty \ep^{\theta}\tau_i\chi_i^\ep(x).\end{array}\]
This choice is made so that for all $y,$ $k^\ep (\cdot,y)$ is a probability measure on $[0,y].$

\subsection{Compactness assumptions on the coefficients}
%%%%%%%%%%%%%%%%%%%%%%%%%%%%%%%%%%%%%%%%%%%%%%%%%%%%%%%%%%%%%%%%%%%%%%%%%%%%%%%%%%%%%%%%%%%%%%%%%%%%%%%%%%%%%%%%%%%%%%%%%%%%%%%%%%%%%%%%%%%%%%%%%%%%%%%%%%%%%%%%%%%%%%%%%%%

For technical purposes
we need further assumptions on the discrete coefficients.
Let us collect them as follows:
\beq\label{as:coef}
\left\{\begin{array}{l}
\textrm{There exists $K>0$ such that }
\\
\big|\beta_{i+1}-\beta_i\big|\leq {Ki^{\alpha-1}}\\
\big|\mu_{i+1}-\mu_i\big|\leq Ki^{m-1},\\
 \big|\tau_{i+1}-\tau_i\big|\leq Ki^{\theta-1},
\end{array}\right.\eeq
where the exponents $\alpha, \theta, m$ are defined in \eqref{growthcond}.
For the fragmentation kernel we assume furthermore
\beq\label{as:k_compact}
\left\{\begin{array}{l}
\textrm{There exists $K>0$ such that for any $i,j$}
\\
 \Big|\dis\sum_{p=0}^{i-1}\dis\sum_{r=0}^{p-1}k_{r,j+1}-\dis\sum_{p=0}^{i-1}\dis\sum_{r=0}^{p-1}k_{r,j}\Big|\leq K.
\end{array}\right.
\eeq
These assumptions 
will be helpful for investigating the behavior of \eqref{eq:discrete_epsilon} as $\ep$ goes to 0 since they provide compactness properties.
We summarize these properties  in the following lemmata.

\begin{lemma}\label{lem:compacite1}
Let $\big(z_i\big)_{i\in \mathbb N}$ be a sequence of nonnegative real numbers verifying
\[
0\leq z_i\leq Ki^\kappa,\qquad
\big|z_{i+1}-z_i\big|\leq {Ki^{\kappa-1}}
\]
for some $K>0$ and $\kappa\geq 0$.
For $x\geq 0$, we set $z^\ep(x)=\sum_{i} \ep^\kappa z_i\chi_{[\ep i,\ep(i+1))}(x)$.
Then there exist a subsequence $\ep_n\rightarrow 0$, and a continuous function $z:x\in [0,\infty)\mapsto z(x)$ such that  
$z^{\ep_n}$ converges to $z$ 
uniformly on $[r,R]$ 
for any $0<r<R<\infty$.
If $\kappa>0$, the convergence holds 
 on $[0,R]$ for any $0<R<\infty$ and we have $z(0)=0$.
\end{lemma}

%MD CHANGE 06 DEC 08
We shall apply this statement to the sequences $\beta^\ep,\mu^\ep,\tau^\ep$.
% \ref{app:compact}.
A similar compactness property can be obtained for the fragmentation coefficients.

\begin{lemma}\label{lem:compacitek}
Let the coefficients $k_{i,j}$ satisfy Assumptions \fer{as:k_symmetric},\fer{as:k_proba}  and \fer{as:k_compact}. Then there exist a subsequence $\big(\ep_n\big)_{n\in \mathbb N}$ and $k:y\in [0,\infty)\mapsto k(\cdot,y)\in
{\cal M}^1_+([0,\infty))$ which belongs to
$ {\cal C}([0,\infty);{\cal M}^1_+([0,\infty))-\mathrm{weak}-\star)$  satisfying also  \fer{as:k_symmetric} and \fer{as:k_proba} (in their continuous version) and such that $k^{\ep_n}$ converges to $k$ in the following sense:
for every compactly supported smooth function $\varphi \in {\cal C}_c^\infty([x_0,\infty))$, denoting 
\beq \label{def:phi}
\phi^{\ep_n}(y)=\int_{n_0(\ep_n)\ep_n}^yk^{\ep_n}(x,y)\varphi(x)\ud x,\qquad
\phi(y)=\int_{x_0}^yk(x,y)\varphi(x)\ud x,
\eeq 
 we have $\phi^{\ep_n}\rightarrow \phi$ uniformly locally in $[x_0,+\infty)$.
\end{lemma}

The detailed  proofs of  Lemma \ref{lem:compacite1} and Lemma \ref{lem:compacitek}  are postponed to  Appendix \ref{app:k}.

\subsection{Main results}
We are now ready to state the main results of this article.

\begin{theorem}\label{th:limit}
Assume \eqref{growthcond} and \eqref{as:coef}. Suppose the fragmentation coefficient fulfill \eqref{as:k_struc}--\eqref{as:k_proba}
and
\eqref{as:k_compact}.Then, there exist a subsequence, denoted $\big(\ep_n\big)_{n\in \mathbb N}$, continuous functions $\mu,\tau,\beta$, and a nonnegative measure-valued function $k$ verifying \fer{as:k_symmetric} and \fer{as:k_proba}, such that 
$$\mu^{\ep_n},\tau^{\ep_n},\beta^{\ep_n},k^{\ep_n}\rightarrow \mu,\tau,\beta,k$$
in the sense of Lemma \ref{lem:compacite1} and Lemma \ref{lem:compacitek}.

Let the initial data satisfy \eqref{cieps}.
Then we can choose the subsequence $\big(\ep_n\big)_{n\in \mathbb N}$ such that there exists $(U,V)$  for which
$$\left\lbrace
\begin{array}{l}
u^{\ep_n} \rightharpoonup U, 
\text{ in ${{\cal C}} ([0,T];{{\cal M}}^1 ([0,\infty))-\mathrm{weak}-\star))$}, \\ 
v^{\ep_n} \rightharpoonup V  \quad \text{ uniformly on $[0,T]$}. 
\end{array}\right.$$
We have $xU(t,x)\in {{\cal M}}^1 ([0,\infty))$, 
the measure $U(t,.)$ has its support included in $[x_0,+\infty)$ for all time $t\geq 0$, 
  and  $(U,V)$ satisfies \eqref{eq:PDEweak}--\eqref{eq:mass_continu}.
%is a monomer preserving weak solution of the prion proliferation equations.
\end{theorem}

\begin{theorem}\label{th:ODE}
The limit $(U,V)$ exhibited in   Theorem~\ref{th:limit}  is a monomer preserving weak solution (\emph{i.e.} satisfies also Equation \fer{eq:ODE}) in the following situations:
\begin{itemize}
\item[i)] $x_0=0$ and either $\theta>0$ (so that  the limit $\tau$ satisfies $\tau(0)=0$),
or the rates $\tau_i=\tau$ are constant.
\item[ii)] $x_0>0$ and the discrete fragmentation coefficients fulfill  the following strengthened   assumption: 
for any $i,j$ we have

\beq \label{as:kij:strong}
 \bigg|\sum_{i'\leq i}\Big(k_{i'j+1}-k_{i',j}\Big)\bigg|\leq \frac{K}{j},
 \qquad k_{i,j}\leq \frac{K}{j}.
\eeq 
\end{itemize}
\end{theorem}       

 %%%%%%%%%%%%%%%%%%%%%%%%%%%%%%%%%%%%%%%%%%%%%%%%%%%%%%%%%%%%%%%%%%%%%%%%%%%%%%%%%%%%%%%%%%%%%%%%%%%%%%%%%%%%%%%%%%%%%%%%%%%%%%%%%%%%%%%%%%%%%%%%%%%%%%%%%%%%%%%%%%%%%%%%%%%%

%%%%%%%%%%%%%%%%%%%%%%%%%%%%%%%%%%%%%%%%%%%%%%%%%%%%%%%%%%%%%%%%%%%%%%%%%%%%%%%%%%%%%%%%%%%%%%%%%%%%%%%%%%%%%%%%%%%%%%%
%
%                                          SECTION:MOMENTS ESTIMATES
%
%
%%%%%%%%%%%%%%%%%%%%%%%%%%%%%%%%%%%%%%%%%%%%%%%%%%%%%%%%%%%%%%%%%%%%%%%%%%%%%%%%%%%%%%%%%%%%%%%%%%%%%%%%%%%%%%%%%%

\section{Moment estimates}\label{sec:moments} 

We start by establishing
a priori estimates 
uniformly with respect to $\ep$.
These estimates will induce compactness properties on the sequence of solutions. As described in \cite{Laurencot-multiple} for general coagulation fragmentation models, the model has the property of propagating moments. 

%%%%%%%%%%%%%%%%%LEMMME%%%%MOMENTS ESTIMATES%%%%%%%%%%%%%%%%%%%%%%%%%%%%%%

\begin{lemma}\label{lem:moment}
Let the assumptions
of Theorem \ref{th:limit} be fulfilled.
 Then for any $T>0$, there exists a constant $C<\infty$ which only depends on $M_0, M_{1+\sigma}, K$ and $T$, such that for any $\varepsilon>0$:
$$\dis\sup_{t\in [0,T]} \int_0^\infty (1+x+x^{1+\sigma})u^{\varepsilon}(t,x)\ud x\leq C,
\qquad 
0\leq v^\ep(t)\leq C$$
%and furthermore
%\[\dis\sup_{t\in [0,T]}
%\dis\int_0^t\dis \int_0^\infty 
% \mu^\ep(x) x^{1+\sigma}u^\ep(s,x)\ud x\ud s.
%\]
\end{lemma}

%%%%%%%%%%%%%%%%%%%%%%%%%%%%%%%%%%%%%%%%%%%%%%%%%%%%%%%%%%%%%%%%%%%%%%%%%%

%%%%%%%%%%%%%%%%%PREUVE%%%%MOMENTS ESTIMATES%%%%%%%%%%%%%%%%%%%%%%%%%%%%%%

%%%%%%%%%%%%%%%%%%%%%%%%%%%%%%%%%%%%%%%%%%%%%%%%%%%%%%%%%%%%%%%%%%%%%%%%%%
\proof 
For $r \geq 0$, we denote 
$$M^\varepsilon_r (t)= \varepsilon\dis\sum_{i=n_0}^\infty (i\varepsilon)^r\ u_i^\varepsilon (t).$$
As in \cite{Thierry}, we can notice that 
$$\int_0^\infty \Big(\frac{x}{2}\Big)^r u^\varepsilon(t,x)\ud x \leq M^\varepsilon_r (t)\leq \int_0^\infty x^r u^\varepsilon(t,x)\ud x.$$
Therefore, we only need to control $M^\varepsilon_0(t), \; M^\ep_1 (t)$ and $M^\varepsilon_{1+\sigma}(t).$
We notice the obvious but useful inequality, for $0\leq r\leq 1+\sigma,$
$$(i\ep)^r\leq 1+(i\ep)^{1+\sigma},$$
and therefore,
$$|M^\ep_r|\leq |M^\ep_0|+|M^\ep_{1+\sigma}|.$$
In the sequel, we use alternatively two equivalent discrete weak formulations of Equation (\ref{eq:discrete_epsilon}) in the spirit of \cite{Well-posed}.  We multiply the second equation of (\ref{eq:discrete_epsilon}) by $\varphi_i$ and summing over $i$, we first obtain
\begin{equation} \begin{array}{lll}
\dfrac{\ud}{\ud t} \dis\sum_{i=n_0}^\infty u_i^\ep\varphi_i&=&-\ep^m\dis\sum_{i=n_0}^\infty\mu_i u_i^\ep\varphi_i-\ep^\alpha\dis\sum_{i=n_0}^\infty\beta_iu_i^\ep\varphi_i\\
&&-\ep^{\theta-1} \dis\sum_{i=n_0}^\infty v^\ep(\tau_iu_i^\ep-\tau_{i-1}u_{i-1}^\ep)\varphi_i+2\ep^\alpha\dis\sum_{i=n_0}^\infty\varphi_i\sum_{j>i}\beta_jk_{i,j}u_j^\ep,\\
&=& -\ep^m\dis\sum_{i=n_0}^\infty\mu_i u_i^\ep\varphi_i-\ep^\alpha\dis\sum_{i=n_0}^\infty\beta_i u_i^\ep\varphi_i+
\ep^{\theta-1} \dis\sum_{i=n_0}^\infty\tau_iu_i^\ep(\varphi_{i+1}-\varphi_i)
\\
&&
+2\ep^\alpha\dis\sum_{i=n_0}^\infty\varphi_i\dis\sum_{j>i}\beta_jk_{i,j}u_j^\ep.\label{weak:1}
\end{array}\end{equation}
Using the properties of $k_{i,j}$, we rewrite  the fragmentation terms as follows
\begin{eqnarray*}
\dis\sum_{i=n_0}^\infty \beta_i u^\ep_i\varphi_i&=&2\dis\sum_{j=n_0+1}^\infty \beta_j\dis\sum_{i=1}^{j-1}ik_{i,j} u^\ep_j\frac{\varphi_j}{j}+\beta_{n_0}u^\ep_{n_0}\varphi_{n_0}\\
																							 &=&2\dis\sum_{j=n_0+1}^\infty \dis\sum_{i=n_0}^{j-1}ik_{i,j}\beta_j u^\ep_j\frac{\varphi_j}{j}+2\dis\sum_{j=n_0+1}^\infty \dis\sum_{i=1}^{n_0-1}ik_{i,j}\beta_j u^\ep_j\frac{\varphi_j}{j}+\beta_{n_0}u^\ep_{n_0}\varphi_{n_0},\\
2\dis\sum_{i=n_0}^\infty\varphi_i\sum_{j>i}\beta_jk_{i,j}u^\ep_j&=& 	2\dis\sum_{j=n_0+1}^\infty\sum_{i=n_0}^{j-1}ik_{i,j}\beta_ju^\ep_j\frac{\varphi_i}{i}.
\end{eqnarray*}
By using \eqref{as:k_mass}, we obtain
\begin{eqnarray*}
2\dis\sum_{i=n_0}^\infty\varphi_i\sum_{j>i}\beta_jk_{i,j}u^\ep_j-\dis\sum_{i=n_0}^\infty \beta_i u^\ep_i\varphi_i&=&-2\dis\sum_{j=n_0}^\infty \dis\sum_{i=1}^{n_0-1}ik_{i,j}\beta_j u^\ep_j\frac{\varphi_j}{j}\\
&&+2\dis\sum_{j=n_0+1}^\infty\dis\sum_{i=n_0}^{j-1}ik_{i,j}\beta_ju^\ep_j\bigg(\frac{\varphi_i}{i}-\frac{\varphi_j}{j}\bigg).
\end{eqnarray*}
Replacing in the weak formulation we get
\begin{equation}\begin{array}{lll}\label{eq:weak2}
\dfrac{\ud}{\ud t} \dis\sum_{i=n_0}^\infty u^\ep_i\varphi_i&=& -\ep^m\dis\sum_{i=n_0}^\infty\mu_i u^\ep_i\varphi_i+
\ep^{\theta -1} v^\ep \dis\sum_{i=n_0}^\infty\tau_iu^\ep_i(\varphi_{i+1}-\varphi_i)
\\
&&+2\ep^\alpha\dis\sum_{j=n_0+1}^\infty\dis\sum_{i=n_0}^{j-1}ik_{i,j}\beta_ju^\ep_j\bigg(\frac{\varphi_i}{i}-\frac{\varphi_j}{j}\bigg)-2\ep^\alpha\dis\sum_{j=n_0}^\infty \dis\sum_{i=1}^{n_0-1}ik_{i,j}\beta_j u^\ep_j\frac{\varphi_j}{j}.
\end{array}
\end{equation}
This last formulation makes the estimates straightforward (the computations are formal but can be understood as uniform bounds on solutions of truncated systems and therefore on any admissible solution).
Taking $\phi_i=i\ep,$ we obtain the first moment, that is, the previously seen mass balance:
\beq \dfrac{\ud}{\ud t}\bigg(v^\ep+\varepsilon^2\dis\sum_{i=n_0}^\infty i u^\ep_i\bigg)=-\gamma v^\ep-\varepsilon^{2+m}\dis\sum_{i=n_0}^\infty \mu_i iu^\ep_i+\lambda\leq\lambda.\label{eq:mass_discret_epsilon}\eeq
Therefore, we get ($u_i^\ep$ and $v^\ep$ are nonnegative)
$$
0\leq v^\ep(t)+M^\varepsilon_1(t)\leq \rho^0+\lambda T \qquad \text {for } 0\leq t\leq T<\infty
$$
and 
\[
\dis\int_0^t \varepsilon^{2+m}\dis\sum_{i=n_0}^\infty \mu_i iu^\ep_i(s,x)\ud s
\leq \rho^0+\lambda T \qquad \text {for } 0\leq t\leq T<\infty.
\]
To obtain an estimate on the $0$th order moment, we take $\varphi_i=\varepsilon$. The term with $\tau_i$ vanishes. Considering only the nonnegative part of the derivative, we derive  from  \fer{eq:weak2}
\begin{eqnarray*}
\frac{\ud}{\ud t}M^\varepsilon_0(t)  &\leq& 2\varepsilon^{1+\alpha} \dis\sum_{j=n_0+1}^\infty\sum_{i=n_0}^{j-1}ik_{i,j}\beta_ju^\ep_j\frac{1}{i},\\
																										&	\leq & 2\varepsilon^{1+\alpha} 	\dis\sum_{j=n_0+1}^\infty\beta_ju^\ep_j\leq 2KM^\varepsilon_\alpha(t).
																										\end{eqnarray*} 
To give the bound on the $(1+\sigma)$th moment, we choose $\varphi_i=\varepsilon (\varepsilon i)^{1+\sigma}$ in the weak formulation. Thanks to the mean value inequality, we have
$$((\varepsilon (i+1))^{1+\sigma}-(\varepsilon i)^{1+\sigma})\leq (1+\sigma)\ep(\ep(i+1))^\sigma
\leq (1+\sigma)2^\sigma\ep (\ep i)^\sigma,$$
therefore  \fer{eq:weak2} yields
\begin{eqnarray*}
\frac{\ud}{\ud t} M^\varepsilon_{1+\sigma}(t) 
+\ep^{1+m}
\dis\sum _{i=n_0}^\infty \mu_i (\ep i)^{1+\sigma}u^\ep_i
&\leq& v^\ep(1+\sigma)2^\sigma\dis\sum_{i=n_0}^\infty \ep^\theta\tau_iu^\ep_i\ep (\ep i)^\sigma,
\\
&\leq&
 K(\rho^0+\lambda T)(1+\sigma)2^\sigma M^\varepsilon_{\theta+\sigma}(t).
\end{eqnarray*}
Since $0\leq \theta\leq 1$, denoting $C=\max(K(\rho^0+\lambda T)(1+\sigma)2^\sigma,2K)$, it leads to
$$\dfrac{\ud}{\ud t}\bigg(M_0^\ep(t)+M^\ep_{1+\sigma}(t)\bigg)\leq C\bigg(M_\alpha^\ep(t)+M^\ep_{\theta+s}(t)\bigg)\leq 2C\bigg(M_0^\ep(t)+M^\ep_{1+\sigma}(t)\bigg),$$
and we conclude by the Gronwall lemma.
It ends the proof of Lemma \ref{lem:moment}. 
\qed

Hereafter, we denote by $C$ a constant depending only on $T, M_0,\rho^0, M_{1+\sigma},K$ and $\lambda$ such that 
$$M_0^\varepsilon,v^\ep,M_1^\varepsilon,M_{1+\sigma}^\varepsilon\leq C.$$
%%%%%%%%%%%%%%%%%%%%%%%%%%%%%%%%%%%%%%%%%%%%%%%%%%%%%%%%%%%%%%%%%%%%%%%%%%

%%%%%%%%%%%%%%%%%LEMMME%%%%v is equicontinuous%%%%%%%%%%%%%%%%%%%%%%%%%%%%

%%%%%%%%%%%%%%%%%%%%%%%%%%%%%%%%%%%%%%%%%%%%%%%%%%%%%%%%%%%%%%%%%%%%%%%%%%
\begin{lemma}\label{lem:equicontinuity}
Under the assumptions of Lemma \ref{lem:moment}, the sequence of monomers concentration $(v^\varepsilon)_{\varepsilon>0}$ is equicontinuous on $[0,T]$.
\end{lemma}
%%%%%%%%%%%%%%%%%%%%%%%%%%%%%%%%%%%%%%%%%%%%%%%%%%%%%%%%%%%%%%%%%%%%%%%%%%

%%%%%%%%%%%%%%%%%Preuve%%%%v is equicontinuous%%%%%%%%%%%%%%%%%%%%%%%%%%%%

%%%%%%%%%%%%%%%%%%%%%%%%%%%%%%%%%%%%%%%%%%%%%%%%%%%%%%%%%%%%%%%%%%%%%%%%%%
\proof We use the estimates of Lemma \ref{lem:moment} to evaluate the derivative of $v^\ep$. We recall the equation satisfied by $v^\ep$
$$\dfrac{\ud v^\ep}{\ud t}=\lambda-\gamma v^\ep +\varepsilon^{1+\theta} v^\ep\sum\tau_i u^\ep_i +2\varepsilon^{2+\alpha}\sum_{i\geq n_0}\sum_{j<n_0}jk_{j,i}\beta_iu^\ep_i,$$
which implies
$$
\left|\dfrac{\ud v^\ep}{\ud t}\right|\leq \lambda+\gamma C+KC^2 +2\varepsilon n_0(\ep)\ KM_\alpha^\varepsilon.$$
Since the sequence $\big(M^\ep_\alpha\big)_{\ep>0}$ is uniformly bounded with respect to $\ep$ by Lemma \ref{lem:moment} (recall that $\alpha \leq 1 + \sigma$), the sequence $\big(v^\ep\big)_{\ep>0}$ 
satisfies a uniform  Lipschitz criterion on  $[0,T]$. 
This concludes the proof of Lemma \ref{lem:equicontinuity}.
\qed

\

{\bf Proof of Theorem~\ref{th:limit}.}
By the Arzela-Ascoli theorem and Lemma \ref{lem:equicontinuity}, there exists a function $V\in {{\cal C}} ([0,T])$ and a subsequence that we still denote $v^{\varepsilon}$ such that
$$v^\ep(t) \longrightarrow V(t)\quad\text{in} \quad {{\cal C}} ([0,T]).$$
In the same way, the moment estimates of Lemma \ref{lem:moment} give uniform boundedness for $(1+x+x^{1+\sigma}) u^{\ep}$ in ${{\cal M}}^1 ([0,\infty)).$ 
Pick a function $\varphi\in {\cal C}^\infty_c([0,\infty))$.  We define
$$\varphi^\ep_i=\int_{i\ep}^{(i+1)\ep}\varphi(x)\ud x,$$
so that $\dis\sum_{n_0\ep}^\infty u^\ep_i\varphi^\ep_i=\int_{0}^\infty u^\ep(t,x)\varphi(x)\ud x,$
and also for $y\in [j\ep,(j+1)\ep[$,
$$\int_0^y k^\ep(x,y)\varphi(x)dx=\int_0^{j\ep}k^\ep(x,j\ep)\varphi(x)dx=\sum_{i=0}^jk_{i,j}\f{\varphi^\ep_i}{\ep}$$
Thanks to the moment estimates of Lemma~\ref{lem:moment}, and using  \fer{weak:1}, we have 
  $$\bigg|\frac{\ud}{\ud t}\int u^{\ep}(t,x)\varphi(x)\ud x\bigg|\leq C(\|\varphi\|_\infty+\|\varphi'\|_\infty)\qquad\text{and}\quad \bigg|\int u^{\ep}(t,x)\varphi(x)\ud x\bigg|\leq C\|\varphi\|_\infty$$
  for some constant $C$ depending only on $K,M_0, M_{1+\sigma},\lambda,T$.
Therefore, for any function $\varphi \in {{\cal C}}_c^\infty ([0,\infty)),$ the integral $\int u^{\ep}(\cdot,x)\varphi(x) \ud x$ is equibounded and equicontinuous. 
   Using a density argument, we can extend this property to $\varphi\in{{\cal C}}_0([0,\infty))$, the space of continuous functions
   on $[0,\infty)$ that tend to 0 at infinity. This means that $\bigl(\int_0^\infty u^\ep(.,x)\varphi(x)\ud x\bigr)_\ep$ belongs to a compact set of ${{\cal C}}(0,T)$. As in \cite{Thierry}, by using the separability of ${{\cal C}}_0([0,\infty))$ and the Cantor diagonal process, we can extract a subsequence $u^{\ep_n}$ and $U\in {{\cal C}}([0,T];{{\cal M}}^1([0,\infty))-\mathrm{weak}-\star)$, such that the following convergence
$$\int_0^\infty u^{\ep_n}(t,x)\varphi(x)\ud x\rightarrow \int_0^\infty U(t,x)\varphi(x)\ud x,$$
as $\ep_n\rightarrow 0$, holds uniformly on $[0,T]$, for any $\varphi\in{{\cal C}}_0([0,\infty))$. 
As $u^\ep(t,x)=0$ for $x\leq \ep n_0(\ep)$, we check that $U(t,.)$ has its support in $[x_0,\infty[$. 
It remains to prove that $(U,V)$ satisfies \eqref{eq:PDEweak} \eqref{eq:mass_continu}.%is a monomer preserving weak solution of the prion proliferation equations.

Let $\varphi$ be a smooth function supported in $[\delta,M]$ with $x_0<\delta<M<+\infty$,  choosing $\ep n_0(\ep)+2\ep <{\delta}$ (what is possible due to  (\ref{as:n0})).
By using Lemma~\ref{lem:compacite1} and Lemma~\ref{cvmt}, we check that, for a suitable subsequence, one has
%
%Moreover, since $(1+x^{1+\sigma})u^\ep$
%is bounded in $L^\infty(0,T;{{\cal M}}^1(0,\infty))$, the same convergence holds for $ x^ru^\ep$ for any $r\in[0,1+\sigma)$. That is 
%\beq\label{eq:xr}
%\forall \varphi\in{{\cal C}}^0([0,\infty)),\forall r\in [0,1+s[ \qquad \int x^ru^{\ep_n}(t,x)\varphi(x)\ud x\rightarrow \int x^rU(t,x)\varphi(x)\ud x\quad \text{in} \; {{\cal C}}([0,T]).
%\eeq
%Using (\ref{eq:xr}) and Lemma~\ref{lem:compacite1}, as $m,\alpha,\theta\in[0,1+\sigma[,$ up to a subsequence, one has for any $\varphi\in{{\cal C}}^c((0,\infty))$ 
%(with compact support in, say, $[r,R]$, $0<r<R<\infty$)
\beq\label{as:convbeta}
\begin{array}{lcl}
\dis\int_0^\infty \mu^{\ep_n}(x)u^{\ep_n}(t,x)\varphi(x)\ud x&\xrightarrow[\ep_n\rightarrow 0]{} &\dis\int_0^\infty \mu(x)U(t,x)\varphi(x)\ud x,\\
[0.3cm]
\dis\int_0^\infty \beta^{\ep_n}(x)u^{\ep_n}(t,x)\varphi(x)\ud x&\xrightarrow[\ep_n\rightarrow 0]{} &\dis\int_0^\infty \beta(x)U(t,x)\varphi(x)\ud x,
\\
[0.3cm]
\dis\int_0^\infty \tau^{\ep_n}(x)u^{\ep_n}(t,x)\varphi(x)\ud x&\xrightarrow[\ep_n\rightarrow 0]{} &\dis\int_0^\infty \tau(x)U(t,x)\varphi(x)\ud x,
\end{array}\eeq
uniformly on $[0,T]$.
Equation (\ref{weak:1}) can be recast in the following integral form
\beq\begin{array}{ll}\label{weak:int}
\dis\f{\ud}{\ud t}\int_{0}^\infty u^\ep (t,x) \varphi (x) \ud x = - \int _{x_0}^\infty \mu^\ep u^\ep (t,x) \varphi (x) \ud x 
-\dis\int_{0}^\infty \tau^\ep u^\ep \Delta^\ep \varphi (x) \ud x \\ \\
- \dis\int_{0}^\infty \beta^\ep u^\ep (t,x) \varphi (x) \ud x 
+ 2 \dis\int_{0}^\infty \int_x^\infty \varphi (x) \beta^\ep(y) u^\ep(t,y) k^\ep (x,y) \ud x \ud y 
\end{array}
\eeq
where we have defined
$$\Delta^\ep \varphi (x) = \int_{i\ep}^{(i+1)\ep}\f{\varphi(s+\ep) - \varphi(s)}{\ep} \;\ud s,\qquad \text{for } x\in[i\ep,(i+1)\ep[.$$

The first and third terms are  treated in \eqref{as:convbeta}. 
Using (\ref{as:convbeta}) again and  remarking that  $|\Delta^\ep(x)-\varphi'(x)|\leq \ep \|\varphi''\|_\infty$, we have 
\beq\label{convtau}
\dis\int_0^\infty \tau^{\ep_n}(x)u^{\ep_n}(t,x)\Delta^{\ep_n}\varphi(x)\ud x\xrightarrow[\ep_n\rightarrow 0]{} \dis\int_0^\infty \tau(x)U(t,x)\varphi'(x)\ud x,
\eeq
uniformly on $[0,T]$.
Let us now study  the convergence of the last term in  \fer{weak:int}.
To this end, we use the notation $\phi$ and $\phi^\ep$ as defined in \fer{def:phi} of Lemma~\ref{lem:compacitek} and we
rewrite 
$$2\int_{x_0}^\infty\int_{x_0}^y\varphi(x)k^\ep(x,y)u^\ep(t,y)\beta^\ep(y)\ud x\ud y=2\int_{x_0}^\infty u^\ep(t,y)\beta^\ep (y)\phi^\ep(y)\ud y.$$
Owing to  (\ref{as:k_compact}) we use Lemma~\ref{lem:compacitek} which leads to
$$\phi^{\ep_n}\xrightarrow[\ep_n\rightarrow 0]{} \phi \qquad\text{uniformly on} \quad [x_0,R] \;\;\text{for any}\;\; R>0,$$ and thus also
$$ \beta^{\ep_n}\phi^{\ep_n}\xrightarrow[\ep_n\rightarrow 0]{} \beta\phi \qquad\text{uniformly on} \quad [x_0,R] \;\;\text{for any}\;\; R>0,$$
for a suitable subsequence.
Finally, we observe  that $\phi^{\ep_n}$ and therefore  $\phi$ are bounded by $\|\varphi\|_{\infty}$. 
Thus, by using the boundedness of the higher order moments of $u^\ep$ in Lemma \ref{lem:moment}
with $1+\sigma> \alpha$, we show  that the fragmentation term passes to the limit (see Lemma \ref{cvmt} in the Appendix).
We finally arrive at 
\beq
\begin{array}{l}
\dis\int_{x_0}^\infty U(t,x)\varphi(x)\ud x-\int_{x_0}^\infty U(0,x)\varphi(x)\ud x\\ \\
= -\dis\int_0^t \int_{x_0}^\infty \mu U (t,x) \varphi (x) \ud x 
-\dis\int_0^tV(s)\int_{x_0}^\infty \tau(x) U(s,x) \varphi' (x) \ud x \\ \\
\qquad- \dis\int_0^t\int_{x_0}^\infty \beta(x) U(s,x) (t,x) \varphi (x) \ud x 
+ 2 \dis\int_0^t\int_{x_0}^\infty \beta(y) U(s,y)\int_0^y\varphi (x)  k (x,y) \ud x \ud y,
\end{array}
\eeq
which is the weak formulation \fer{eq:PDEweak}. Moreover, \fer{eq:mass_discret_epsilon} recasts as 
\[\begin{array}{lll}
v^\ep(t)+\dis\int_0^\infty e^\ep(x) u^\ep(t,x)\ud x &=& v^{0,\ep}(t)+\dis\int_0^\infty e^\ep(x) u^{\ep}(0,x)\ud x
\\
&&+\lambda t - \gamma\dis\int_0^t v^\ep(s)\ud s- 
\dis\int_0^t\dis\int_0^\infty e^\ep(x) \mu^\ep(x) u^\ep(s,x)\ud x\ud s
\end{array}\]
where \[e^\ep(x)=\dis\sum_{i=0}^\infty\ep i\ \chi_{[i\ep,(i+1)\ep)}(x).\]
Clearly $e^\ep(x)$ converges to $x$ uniformly.
Using  the moment estimate in Lemma \ref{lem:moment}, with $\sigma>0$,
we obtain
\[
v^{\ep_n}(t)+\dis\int_0^\infty e^{\ep_n}(x) u^{\ep_n}(t,x)\ud x
\xrightarrow[\ep_n\rightarrow 0]{} V(t)+\dis\int_0^\infty xU(t,x)\ud x
\]
uniformly on $[0,T]$ as well as
\[\dis\int_0^t\dis\int_0^\infty e^{\ep_n}(x) \mu^{\ep_n}(x) u^{\ep_n}(s,x)\ud x\ud s
\xrightarrow[\ep_n\rightarrow 0]{} 
\dis\int_0^t\dis\int_0^\infty x\mu(x) U(s,x)\ud x\ud s.
\]
(We refer again to  Lemma \ref{cvmt}, or a slight adaptation of it.) 
As
 $\ep_n\rightarrow 0$
we are thus led to \eqref{eq:mass_continu}.
\qed

\

{\bf Proof of Theorem~\ref{th:ODE}.}
We rewrite the rescaled ODE as 
$$\f{d v^\ep}{dt}=\lambda-\gamma v^\ep -\int_{n_0\ep}^\infty \tau^\ep(x)u^\ep(t,x)dx+2\int_{n_0\ep}^\infty \beta^\ep(y)u^\ep(t,y)\int_0^{n_0^\ep}e^\ep(x)k^\ep(x,y)dx,$$
Depending on the value of $x_0$, we have to care about the last term ($x_0>0$) or the next two last term ($x_0=0$).
As already remarked in the proof of Lemma \ref{lem:equicontinuity}, in case where $x_0=0$,
the fragmentation term can be dominated by
\[
2\varepsilon^{2+\alpha}\dis\sum_{i\geq n_0}\sum_{j<n_0}jk_{j,i}\beta_iu^\ep
\leq 2\varepsilon n_0(\ep)\ KM_\alpha^\varepsilon.\]
Hence this contribution vanishes as $\ep$ goes to 0 when $\lim_{\ep \rightarrow 0}\ep n_0(\ep)=x_0=0$.
 Nevertheless for case i) we still have to justify that  $\int_0^\infty \tau^\ep(x)u^\ep (t,x)\ud x$
 passes to the limit. We get 
 \beq\label{convtaueqv}
 \int_{n_0(\ep_n)\ep_n}^\infty\tau^{\ep_n}u^{\ep_n}(t,x)\ud x\xrightarrow[\ep_n\rightarrow 0]{} \int_{x_0}^\infty\tau U(t,x)\ud x,\qquad \text{in}\quad {{\cal C}}([0,T])\eeq
 by using the strengthened assumption  $0<\theta\leq 1$ in \eqref{cieps}.
 Indeed it implies that $\tau^\ep(x)$ converges uniformly to $\tau(x)$ on any compact set $[0,R]$ while  these functions do not grow 
faster than $x$ at infinity. We can thus use Lemma \ref{cvmt} to conclude.
\

In the situation ii), 
another difficulty comes from the fragmentation term since we have to prove that
\[
2\int_{n_0(\ep_n)\ep_n}^\infty\int_{0}^{n_0(\ep_n)\ep_n}
e^{\ep_n}(x) k^{\ep_n}(x,y) \beta^{\ep_n}(y)u^{\ep_n}(t,y) \ud x \ud y
\xrightarrow[\ep_n\rightarrow 0]{}
2\int_{x_0}^\infty\int_{0}^{x_0}xk(x,y)\beta(y)U(t,y)\ud x\ud y.\]
 The stronger compactness assumptions  \eqref{as:kij:strong} are basically Ascoli-type assumptions on the repartition function associated to the kernels $k^\ep$. Denoting, in a similar manner to Appendix \ref{app:k}:
$$F^\ep (x,y)= \int_0^x k^\ep (z,y) \ud z,\qquad G^\ep (x,y)=\int_0^x F^\ep (z,y)\ud z,$$
Lemma \ref{lem:th3} (see Appendix \ref{app:k}) ensures that  $F^\ep\rightarrow F$ uniformly  on compact sets of $\R_+\times[x_0,+\infty)$. We also make the remark that 
 $$\Big|\int_0^{n_0\ep}e^\ep(x)k^\ep(x,y)\ud x-\int_0^{n_0\ep}xk^\ep(x,y)\ud x\Big|\leq \ep,$$
 $$\int_0^{n_0\ep}xk^\ep(x,y)\ud x=\Big[xF^\ep(x,y)\Big]^{x=n_0\ep}_{x=0}-\int_0^{n_0\ep}F^\ep(x,y)=(n_0\ep)F^\ep(n_0\ep,y)-G^\ep(n_0\ep,y).$$
 Thanks to Lemma \ref{lem:th3}, we know that the concerned quantities are uniformly bounded and converge uniformly on compact sets, so that
 $$\int_0^{n_0\ep_n}e^{\ep_n}(x)k^{\ep_n}(x,y)\ud x\xrightarrow[\ep_n\rightarrow 0]{}
 \int_0^{x_0}xk(x,y)dx \qquad\text{uniformly on compact sets}.$$
 And as before this is sufficient to prove that 
 $$2\int_{n_0\ep_n}^\infty \beta^{\ep_n}(y)u^{\ep_n}(t,y)\int_0^{n_0 \ep_n}e^{\ep_n}(x)k^{\ep_n}(x,y)\ud x
\ud y
\xrightarrow[\ep_n\rightarrow 0]{}
 2\int_{x_0}^\infty \beta(y)u(t,y)\int_0^{x_0}xk(x,y)\ud x\ud y.$$\qed
 \section{Boundary Condition for the Continuous System}
\label{sec:bound}

The discrete system \fer{eq:discret} only needs an initial condition 
prescribing the $u_i$'s and $v$ at time $t=0$
to be well-posed, as Theorem \ref{th:exist:discret} states. It is different for the continuous system \fer{eq:ODE}\fer{eq:PDEstrong}: a boundary condition at $x=x_0$ is needed when $\tau(x_0)>0$
(in which case the characteristics associated to the ``velocity'' $\tau$ are ``incoming'').
Even when $\tau(x_0)=0$, difficulties might arise when $x\mapsto \tau(x)$ is not regular enough to define the associated characteristics. 
However, according to the analysis of \cite{DubSte, Well-posed}, we have seen 
that 
the notion of  ``monomer preserving  solution'' 
appears naturally, inserting  \eqref{eq:balance} as a constraint.
It leads to the question of deciding how this condition is related to   \eqref{eq:ODE} and
\eqref{eq:PDEstrong} 
and to determine the corresponding boundary condition to be used at $x=x_0$.

Let $(U,V)$ be a   ``monomer preserving'' solution. In this section we do not care about regularity requirement, and we perform several manipulations on the solution (that is assuming all the necessary integrability conditions).
We suppose that the kernel $k$ splits into a Dirac mass at $x=x_0$ and a measure which is diffuse at $x_0$:
\[k(x,y)=l(x,y)+\delta(x=x_0^+)\psi^+(y)+\delta(x=x_0^-)\psi^-(y),\]
where for any $y\geq 0$, $\int_{x_0-\eta}^{x_0+\eta} l(x,y)\ud x\rightarrow 0$ as $\eta$ goes to 0.
We have to distinguish between $x_0^+$ and $x_0^-$ since their biological and mathematical interpretation is different: the Dirac mass at $x_0^+$ means that polymers of size $x_0$ are formed, whereas the Dirac mass at $x_0^-$ is interpreted as breakages of polymers of size $x_0$ going back to the monomers compartment $V.$ As shown below, the mathematical treatment of each is different.
The time derivative of \fer{eq:mass_continu} leads to
\[%\begin{array}{l}
\dis\frac{\ud }{\ud t}\varrho
=
\dis\frac{\ud V}{\ud t}+
\dis\int_{x_0}^\infty x \f{\p}{\p t} U (t,x) \ud x =-
\dis\int_{x_0}^\infty x\mu(x) U(t,x) \ud x
\\
+\lambda -\gamma V
%-v\int_{x_0}^\infty \tau(x) U(t,x) \ud x+ 2\int_{x_0}^\infty \int_0^{x_0} y k(y,x) \beta (x) U(t,x) \ud y \ud x
%\end{array}
\]
In the left hand side, we can compute 
the derivative 
of the moment of $U$ by using \eqref{eq:PDEstrong}. We get
$$\begin{array}{lll}
\dis\frac{\ud }{\ud t}\dis\int_{x_0}^\infty xU(t,x)\ud x
&=&
-\dis\int_{x_0}^\infty x\beta U \ud x 
-\dis\int_{x_0}^\infty x\mu U \ud x
\\
&&
- V\dis\int_{x_0}^\infty x \f{\p}{\p x} (\tau U) \ud x
+2\dis\int_{x_0}^\infty x \int_{x}^\infty l(x,y) \beta(y) U(t,y) \ud y \ud x.
\end{array} $$
In this equation, since \fer{eq:PDEstrong} is only written for $x>x_0,$ only the diffuse part of the kernel $k$ appears.
Integrating by parts, the convection term yields
\[
\dis\int_{x_0}^\infty x \f{\p}{\p x} (\tau U) \ud x=-x_0\tau(x_0)U(t,x_0) - 
\dis\int_{x_0}^\infty \tau U \ud x.
\]
Now we use \eqref{eq:ODE}, which writes
\[ \begin{array}{ll}
\dis\frac{\ud V}{\ud t}=  \lambda-\gamma V -V\dis\int_{x_0}^\infty \tau(x)U(t,x)\ud x &+2\dis\int_{x=x_0}^\infty\int_{y=0}^{x_0}y\ l(y,x)\beta(x)U(t,x)\ud y\ud x
 \\
& + 2 x_0 \dis\int_{x=x_0}^\infty \psi^-(x)\beta(x)U(t,x)\ud x
\end{array}
\]
then we obtain
$$
x_0V(t)\tau(x_0)U(t,x_0)
-\dis\int_{x_0}^\infty x\beta(x)
U(t,x)\ud x
+2\dis\int_{x_0}^\infty \beta(x)U(t,x)\left( \int_{0}^x yl(y,x)  \ud y + x_0 \psi^-(x) \right)\ud x
=0
.
$$
However, \eqref{as:k_mass} is interpreted as
\[2\dis\int_0^x yl(y,x)\ud y + 2 x_0\chi_{[x_0,\infty)}(x)\psi^-(x)+ 2 x_0\chi_{(x_0,\infty)}(x)\psi^+(x)=x.\]
We are thus led to the relation:
\[
x_0\left(V(t) \tau (x_0) U(t,x_0)- 2\dis\int_{x_0}^\infty \psi^+(x) \beta(x)U(t,x) \ud x\right)=0
\]
which suggests the boundary condition
\beq \label{eq:boundary:general}
 x_0 V \tau (x_0) U(t,x_0) = 2 x_0 \int_{x_0}^\infty \psi^+(x) \beta(x) U(t,x) \ud x.
\eeq 
(Note that written in this way it makes sense also when $x_0=0$.)

When $x_0>0$, %no boundary condition can be derived from the previous results in a completely rigorous way, since we are not able to say that Equation \fer{eq:ODE} is verified. However, 
the above calculation gives solid intuitive ground to choose Equation \fer{eq:boundary:general} as a boundary condition, defining the incoming flux by means of a 
weighted average 
of the solution over the size variable.
In particular %for $\tau(x_0)>0$, 
if the Dirac part vanishes we obtain $$V \tau (x_0) U(t,x_0) =0,$$
the boundary condition  proposed in \cite{Greer}, for constant coefficient $\tau$.
It is also the boundary condition used in \cite{Well-posed}.

 %When $x_0=0$ and $\tau(0)=0$,
%as in case i) of Theorem \ref{th:ODE}, there is no necessary boundary condition for Problem \fer{eq:ODE}--\fer{eq:PDEstrong} to be well-posed, and indeed Equations \fer{eq:ODE} and \fer{eq:mass_continu} are formally equivalent, as soon as \fer{eq:PDEstrong} is verified.

If $x_0=0,$ the problem is still harder, since Equation \fer{eq:boundary:general} is empty. Dividing it by $x_0>0$ and passing formally to the limit would however give:
\beq \label{eq:boundary:0}
V \tau (0) U(0) = 2 \int_{0}^\infty \psi^+(x)\beta(x) U(t,x) \ud x.
\eeq 
Here again, it generalizes what has been proposed in \cite{Greer} for $\tau$ constant and $k(x,y)=\f{1}{y} \chi_{x\leq y},$ though without any rigorous justification, and  %, contrarily to the case $x_0>0.$ 
 if $\psi^+=0$ it imposes a vanishing incoming flux.

\section{Discussion on the parameters and choice for $\ep$}
\label{sec:discuss}
\subsection{Orders of magnitude}
A biological discussion upon the parameter values can be found in \cite{Natacha} and is based on \cite{Masel,Masel2005} and references therein.

To carry out the previous scaling limit theorem, we made the following assumptions:
$$s=\frac{\cal U}{\cal V}=\ep^2,\quad \nu=\f{1}{\ep},\quad \displaystyle\lim\ _{\ep\to 0} \ep n_0(\ep)=x_0,\quad \eta=a=c=d=1.$$ 
Let us denote $i_0$ the average size of polymers. Even if there still exists much uncertainty upon its value, we can estimate that the typical size of polymers ranges between $15$ and $1000,$ so we can write  
$$\ep_1 = \frac{1}{i_0} \ll 1.$$
It is also known that the ``conversion rate'' of PrPc is around 5 to 10\% at most (depending on the experiment, on the stage of the disease, etc) ; it means that the mass of proteins present in the monomeric form is much larger than the mass of proteins involved in polymers. In terms of characteristic values, it writes
$$ \ep_2= \frac{i_0 {\cal U}}{\cal V} \ll 1.$$
Finally, we get:
$$\ep = \sqrt{\frac{\cal U}{\cal V}}=\sqrt{\ep_1\ep_2} \ll 1.$$
 Hence, it legitimates the assumption on the parameters $s$ and $\ep.$
Concerning the parameter $a,$ we have $a=\f{L}{\cal V}\approx \f{2400}{500},$ which is in the order of 1.
We have  only $d_0\leq 5.10^{-2}:$ this should lead to neglect the degradation rate of polymers and simplify the equation.

For the fragmentation frequency, it is in the order of the exponential growth rate, found experimentally to be in the order of $0.1;$ in the case of Masel's articles \cite{Masel,Masel2005}, it is supposed that $\alpha=1,$ so it seems relevant (it leads to a fragmentation frequency in the order of $\ep$). However, it has to be precisely compared to the other small parameters which are given by the typical size $i_0$ and the conversion rate to justify the approximation. Moreover, the assumption of a linear fragmentation kernel $\beta$ has to be confrounted to experiments.

Concerning the aggregation rate $\cal T,$ and its related parameter $\nu=\tau V,$ as shown in \cite{Natacha}, in most cases we have $\f{1}{\nu}$ in the range of $[0.01, 0.1],$ so it seems justified to suppose it small ; what has to be explored is its link with the other previously seen small parameters.

To conclude (or open the debate), it seems that each specific experiment, like PMCA protocole, \emph{in vitro} or \emph{in vivo} measures, or yet for the case of recombinant $PrP$ (see \cite{Rezaei08}), the orders of magnitude of each parameter should be carefully estimated, in order to adapt the previous model and stick to the biological reality - which proves to be very different in \emph{in vivo}, \emph{ex vivo} or \emph{in vitro} situations, or yet at the beginning (when there are still very few polymers) and at the end of experiences. The following discussion illustrates this idea, and gives some possible extensions to the previously seen models.

\subsection{Discussion on the fragmentation rates $k_{i,j}$}
To illustrate the central importance of a good estimate of the orders of magnitude, we exhibit here a case where the limit is not the continuous System \fer{eq:ODE}\fer{eq:PDEstrong}, but another one. Our calculation is formal, but a complete proof can be deduced from what preceeds and from \cite{Thierry}.

Let us take, instead of $b=\ep^\alpha:$
$$b=\ep^{\alpha-1},$$
and suppose also that the fragmentation kernel verifies:
$$k_{1,i}=k_{i-1,i}=\f{1}{2}(1-\ep r_i), \quad k_{i,j}=\ep k_{i,j}^0 r_j, 2\leq i \leq j-2.$$
It means that the polymers are much more likely to break at their ends than in the middle of their chain. In that case, under Assumption~\fer{as:k_compact} on $k_{i,j}^0$ and \fer{as:coef} on $r_j$ and $\beta_j,$ if $\alpha -1 \leq 1+\sigma,$

the limit system writes:
\beq 
\left\lbrace\begin{array}{l}\label{eq:continu2}
\dis\frac{\ud v}{\ud t}=\lambda-\gamma v +v\int_{x_0}^\infty \tau(x)U(t,x)\ud x 
\\
\qquad\qquad\qquad
-\dis\int_{x_0}^\infty \beta(x)U(t,x)\ud x +2\int_{x=x_0}^\infty\int_{y=0}^{x_0}yk(y,x)r(x)U(t,x)\ud y\ud x,\\[0.3cm]
\dis\frac{\partial  u}{\partial t}=-\mu(x) U(t,x)-r(x)U(t,x)-v\f{\p}{\p x}(\tau U)+ \f{\p}{\p x} (\beta u) +2\int_x^\infty r(y)k(x,y)U(t,y)\ud y.
\end{array}\right.
\eeq 
Notice also that System~\fer{eq:ODE}\fer{eq:PDEstrong} includes the case of ``renewal'' type equations (refer to \cite{BP} for instance), meaning that the ends of the polymers are more likely to break. For instance, if we have, in the above setting:
$$k_{i-2, i}^0=k_{2, i}^0=\f{m_i}{2},\qquad k_{i, j} ^0 = 0(\f{1}{j}),\quad 3\leq i\leq j-3,$$
then Equation (\ref{eq:continu2}) remains valid, but we have to write the boundary condition as:
\beq \label{eq:bound2}
\tau (x=0) U(x=0)= \int m(y) U(t,y) \ud y.
\eeq 
This indeed can also be written as the measure of $0$ by  $d\mu_y = k(x,y) \ud x:$ $m(y)=\mu_y(\{0\}).$ 

Both of these cases mean that the ends of polymers are more likely to break. What changes is the order of magnitude of what we mean by ``more likely to break'': is it in the order of $\f{1}{\ep},$ in which case System \fer{eq:ODE}\fer{eq:PDEstrong} is valid but with a boundary condition of type (\ref{eq:bound2}) ? Or is the difference of the order of $\f{1}{\ep^2},$ in which case Equation (\ref{eq:continu2}) is more likely ? Refer to \cite{Natacha} for a more complete investigation of what model should be used in what experimental context.

\subsection{Discussion on the minimal size $n_0$}
\label{sec:n0}
We have seen above that to have $x_0=0,$ it suffices to make Assumption (\ref{as:n0}). Having also seen that the typical size $i_0$ is large, and that 
$$\ep^2=\f{1}{i_0} \f{M}{m_1 V},\qquad \f{M}{m_1 V} \ll 1,$$
it is in any case valid to suppose that
$$\f{1}{i_0}=\ep^c,\qquad 0<c<2.$$
Hence, Assumption (\ref{as:n0}) can be reformulated as:
\beq
n_0 \ll i_0^{\f{1}{c}}.
\eeq 
For $c=1,$ it means $n_0 \ll i_0,$ which is true. On the contrary, if we suppose that $x_0>0,$ it means that $n_0 \approx i_0^{\f{1}{c}}:$ in most cases, where for instance $i_0=100$ or $i_0=1000,$ it seems irrelevant. 

%%%%%%%%%%%%%%%%%%%%%%%%%%%%%%%%%%%%%%%%%%%%%%%%%%%%%%%%%%%%%%%%%%%%%%%%%%%%%%%%%%%%%%%%%%%%%%%%%%%%%%%%%%%%%%%%%%%%%%%
%
%                                           APPENDICE
%
%
%%%%%%%%%%%%%%%%%%%%%%%%%%%%%%%%%%%%%%%%%%%%%%%%%%%%%%%%%%%%%%%%%%%%%%%%%%%%%%%%%%%%%%%%%%%%%%%%%%%%%%%%%%%%%%%%%%

\appendix\section{Appendix}

\subsection{Compactness of the coefficients}

\noindent
{\bf Proof of Lemma \ref{lem:compacite1}.} We refer to \cite{Thierry} for the case $\kappa=0$. We prove here the case $\kappa>0$.
First, we show that $z^\varepsilon$ is close to a subsequence satisfying the requirements of the Arzela--Ascoli theorem on $[r,R]$. We define $\tilde z ^\varepsilon$ by 
$$\tilde z ^\varepsilon(x)=\varepsilon^\kappa z _i +\varepsilon^\kappa \frac{z _{i+1}-z _i}{\varepsilon}(x-i\ep) \;\text{for}\; i\ep\leq x\leq (i+1)\ep.$$
We have 
\begin{eqnarray*}
|\tilde z ^\varepsilon(x)- z ^\varepsilon(x)|&=&|\varepsilon^\kappa \frac{z _{i+1}-z _i}{\varepsilon}(x-i\ep)|,\\
&\leq & \varepsilon^\kappa|z _{i+1}-z _i|,\\
&\leq & \varepsilon K(\varepsilon i)^{\kappa-1}\leq 2\varepsilon (Kr^{\kappa-1}+KR^{\kappa-1}).
\end{eqnarray*}
Furthermore $\tilde z ^\varepsilon$ has a bounded derivative since 
\begin{eqnarray*}
\bigg|\frac{\ud \tilde z ^\varepsilon}{\ud x}\bigg|&=&\varepsilon^\kappa \frac{z _{i+1}-z _i}{\varepsilon},\\
&\leq& K(\varepsilon i)^{\kappa -1},\\
&\leq & K r^{\kappa-1}+KR^{\kappa-1}.
\end{eqnarray*}
Therefore, the family $\tilde z ^\varepsilon$ satisfies the requirements of Arzela Ascoli theorem for any interval $[r,R]$ with $0<r<R<+\infty$. We can extract a subsequence converging uniformly to 
$z $. The limit is continuous and satisfies $z (x)\leq K x^\kappa$. When $\kappa>0$ the convergence extends on $[0,R]$ owing to  the remark 
$$\sup_{x\in [0,r]}\big|(z ^\ep-z)(x) \big|\leq 2Kr.$$
This concludes the proof. \qed

During the proof of  Theorem \ref{th:limit} and Theorem \ref{th:ODE} we made repeated
use of the following claim.

\begin{lemma}
\label{cvmt}
Let $z_n$ converge to a continuous function $z$ uniformly on $[0,M]$ for any $0<M<\infty$, with 
$|z_n(x)|\leq K(1+x^\kappa)$.
Let $\big(u_n\big)_{n\in\mathbb N}$ be a sequence of integrable functions which converges to 
$U$ weakly-$\star$ in $\mathcal M^1([0,\infty))$.
We suppose furthermore that 
\[\dis\sup_{n\in\mathbb N}\dis\int_{0}^\infty(1+x^\ell)|u_n(x)|\ud x\leq C<\infty.\]
Assuming $0\leq\kappa<\ell$, we have
\[\dis\int_{0}^\infty z_n(x)u_n(x)\ud x\xrightarrow[n\rightarrow\infty]{}\dis\int_{0}^\infty z(x)U(x)\ud x.\]
\end{lemma}

\proof
Let $\zeta\in C^\infty_c([0,\infty))$ such that $0\leq \zeta(x)\leq 1$,
 $\zeta(x)=1$ on  $[0,R]$, $0< R<\infty$ and $\mathrm{supp}(\zeta)\subset [0,2R]$.
We split 
\[\begin{array}{l}
\Big|
\dis\int_{0}^\infty z_n(x)u_n(x)\ud x-
\dis\int_{0}^\infty z(x)U(x)\ud x
\Big|
\\
=
\Big|\dis\int_{0}^\infty \big(z_n(x)u_n(x)-z(x)U(x)\big)\big(\zeta(x)+1-\zeta(x)\big)\ud x\Big|
\\
\qquad\leq
\dis\int_{0}^\infty |z_n(x)-z(x)|\ |u_n(x)|\ \zeta(x)\ud x
+
\Big| \dis\int_{0}^\infty z(x)\zeta(x) \ \big(u_n(x)-U(x)\big)\ud x\Big|
\\
\qquad\qquad
+\dis\int_{0}^\infty |z_n(x)u_n(x)-z(x)U(x)|\ \big(1-\zeta(x)\big)\ud x.
\end{array}\]
The last integral can be dominated by
\[K\dis\sup_{y\geq R} \left(\dis\frac{1+y^\kappa}{1+y^\ell}\right)\
\left(\dis\sup_{n}\dis\int_{0}^\infty(1+x^\ell)(|u_n(x)|+|U(x)|)\ud x\right).\]
Since $0\leq \kappa<\ell$, this contribution can be made arbitrarily small
by choosing $R$ large enough, uniformly with respect to $n$.
Moreover, we clearly have
\[0\leq \dis\int_{0}^\infty |z_n(x)-z(x)|\ |u_n(x)|\ \zeta(x)\ud x
\leq \dis\sup_{0\leq x\leq 2 R}|z_n(x)-z(x)|\ \dis\sup_{m}\dis\int_{0}^\infty |u_m(x)|\ud x
\xrightarrow[n\rightarrow\infty]{}0
\]
and of course
\[\dis\int_{0}^\infty z(x)\zeta(x) \ \big(u_n(x)-u(x)\big)\ud x
\xrightarrow[n\rightarrow\infty]{}0.\]
Combining all together these informations ends the proof.
\qed

\subsection{Compactness of the fragmentation kernel}
\label{app:k}
We look for conditions on the coefficients guaranteeing some compactness of $k^\ep$. We use a few classical results of convergence of probability measures (see \cite{Billingsley_book} for instance). 
Let us introduce a few notations.
Given a probability-measure-valued function $y\in \mathbb R \mapsto k(.,y)\in {\cal{M}}^1(\mathbb R)$,
we denote $F(.,y)$ its repartition function: $F(x,y)=\int_{-\infty}^x k(s,y)\ud s$ and $G(x,y)$ the function $\int_{-\infty}^x F(z,y)\ud z$. 
We shall deduce the compactness of $k^\ep$ from the compactness of the associated $G^\ep$.
To this end, we need several elementary statements.

\begin{lemma}
Let $\{P^n,\ n\in \mathbb N\}$ be a family of probability measures on $\R$, having their support included in some interval $[a,b]$. We denote $F^n$ the repartition function of $P^n$, and $G^n$ the functions defined by $\int_{-\infty}^x F^n(s)\ud s$. The following assertions are equivalent:
\begin{enumerate}
\item $P^n\rightarrow P$ weakly (\emph{i.e.}, $\forall\;f\in{\cal C}_b (\R),$ $P_n f \rightarrow Pf$)
\item $F^n(x)\rightarrow F(x)$ for all $x$ at which $F$ is continuous
\item $G^n\rightarrow G$ uniformly locally
\end{enumerate}  
\end{lemma}

\begin{lemma}[Conditions for $F$]\label{lem:F}
Let $F$ be a nondecreasing function on $\R$. There exists a unique probability measure $P$ on $\R$, such that $F(x)=P(]-\infty,x])$, iff
\begin{itemize}
\item $F$ is rightcontinuous everywhere,
\item $\lim_{x\rightarrow-\infty} F(x)=0,\lim_{x\rightarrow+\infty} F(x)=1$.
\end{itemize}
Furthermore 
$P$ has its support included in $[a,b]$ iff $F\equiv 0$ on $]-\infty,a[$ and $F(b)=1$. 
\end{lemma}
\begin{lemma}[Conditions for $G$]\label{lem:G}
Let $G$ be a convex function on $\R$. There exists a probability measure $P$ on $\R$, having its support included in $[a,b]$,  such that $G(x)=\int_{-\infty}^x F(s)\ud s$, where $F(x)=P(]-\infty,x])$, iff
\begin{itemize}
\item $G$ is increasing,
\item for $x>b$, $G(x)=G(b)+x$,
\item $G\equiv 0$ on $]-\infty,a]$.
\end{itemize} 
\end{lemma}
\begin{corollary}
Let $\big(G^n\big)_{n\in\mathbb N}$ a sequence satisfying the assumptions of lemma~\ref{lem:G}. Suppose $G^n\rightarrow G$ uniformly locally on $\R$, then $G$ also satisfy these assumptions and we have $P^n\rightarrow P$ weakly.
\end{corollary}
\proof We define the function $F$ as $F(x)=\lim_{\delta \rightarrow 0^+} \f{G(x+\delta)-G(x)}{\delta}$, it is then easy to check that $F$ satisfies assumptions of lemma~\ref{lem:F}, and $G(x)=\int_{-\infty}^x F(s)\ud s$.
\qed
 
%\begin{lemma}\label{lem:compacite_k}
%Suppose that the discrete coefficients satisfy \fer{as:k_compact}.
%Then there exist a subsequence $\ep_n$ and $k\in {\cal {C}}([0,\infty),{\cal {M}}^1_+([0,\infty))-\mathrm{weak}-\star)$ such that 
%\begin{itemize}
%\item $k$ satisfies \fer{as:k_proba},\fer{as:k_symmetric} (and therefore \fer{as:k_mass}),
%\item for every $y>0$, $k^{\ep_n}(.,y)\rightarrow k(.,y)$ in law,
%\item for every $\varphi\in {\cal{C}}^\infty_c([0,\infty))$, $\phi^{\ep_n}\rightarrow \phi$ uniformly on $[0,R]$ for any $0<R<\infty$.
%\end{itemize}
%\end{lemma}

\

{\noindent {\bf Proof of Lemma \ref{lem:compacitek}. } 
We prove the following result, which contains Lemma \ref{lem:compacitek}.

\begin{lemma}\label{lem:compacite_k}
Suppose that the discrete coefficients satisfy \fer{as:k_compact}.
Then there exist a subsequence $\ep_n$ and $k\in {\cal {C}}([0,\infty),{\cal {M}}^1_+([0,\infty))-\mathrm{weak}-\star)$ such that 
\begin{itemize}
\item $k$ satisfies \fer{as:k_proba},\fer{as:k_symmetric} (and therefore \fer{as:k_mass}),
\item for every $y>0$, $k^{\ep_n}(.,y)\rightarrow k(.,y)$ in law,
\item for every $\varphi\in {\cal{C}}^\infty_c([0,\infty))$, $\phi^{\ep_n}\rightarrow \phi$ uniformly on $[0,R]$ for any $0<R<\infty$.
\end{itemize}
\end{lemma}

For any $y\geq 0$, $k^{\ep}(x,y)\ud x$ defines a probability 
measure on $[0,\infty)$, supported in $[0,y]$.
We set 
$F^{\ep}(x,y)=\int_0^xk^{\ep}(z,y)\ud z$ and $G^{\ep}(x,y)=\int_{0}^xF^{\ep}(z,y)\ud z$.
Let $\varphi\in C^\infty_c(\R^*_+)$. 
We start by rewriting, owing to integration by parts, 
$$\phi^{\ep}(y)=\varphi(y)-\int_{0}^y F^{\ep_n}(x,y)\varphi'(x)\ud x =\varphi(y)-G^{\ep}(y,y)\varphi'(y)+\int_{0}^y G^{\ep}(x,y)\varphi''(x)\ud x,$$
where we used the fact that $F^{\ep}(y,y)=\int_{0}^ yk^{\ep}(z,y)\ud z=1$.
The proof is based on the following argument: $G^{\ep}$ is close to a $\tilde G^\ep$ which satisfies the assumptions of the Arzela-Ascoli theorem. 
Given $x,y\geq 0$ and $\ep>0$, $i,j$ denote the integers  satisfying $x\in[i\ep,(i+1)\ep[$, $y\in[j\ep,(j+1)\ep[$
and 
a short computation leads to 
\[\begin{array}{l}F^\ep(x,j\ep)=S_{i,j}+\dis\f{x-i\ep}{\ep}k_{i,j},\\
G^\ep(x,j\ep)=\ep\dis\sum_{p=0}^{i-1}S_{p,j} +(x-i\ep)S_{i,j}+\f{\ep}{2}S_{i,j}+\f{(x-i\ep)^2}{2\ep}k_{i,j},\end{array}\]
where 
\[S_{i,j}=\dis\sum_{r=0}^{i-1}k_{r,j}.\]
We define 
\[\tilde k^\ep(x,y)=\dis\frac{(j+1)\ep-y}{\ep}k^\ep(x,j\ep)+\frac{y-j\ep}{\ep}k^\ep(x,(j+1)\ep)\]
and  we have
$$\tilde G^\ep(x,y)=\frac{(j+1)\ep-y}{\ep}G^\ep(x,j\ep)+\frac{y-j\ep}{\ep}G^\ep(x,(j+1)\ep).$$
Observe 
that 
\[\begin{array}{lll}
|\tilde G^\ep(x,y)-G^\ep(x,y)|&=&
\dis\frac{y-j\ep}{\ep}|G^\ep(x,(j+1)\ep)-G^\ep(x,j\ep)|
\\
&\leq&
\Big|
\ep\dis\sum_{p=0}^{i-1}(S_{p,j+1}-S_{p,j}) +(x-i\ep)(S_{i,j+1}-S_{i,j})+\dis\f{\ep}{2}(S_{i,j+1}-S_{i,j})
\\
&&\qquad\qquad
+\dis\f{(x-i\ep)^2}{2\ep}(k_{i,j+1}-k_{i,j})
\Big|.
\end{array}\]
Due to \eqref{as:k_proba}, we have $0\leq k_{i,j}\leq 1$
and thus $|k_{i,j+1}-k_{i,j}|\leq 1$.
Similarly $ 0\leq S_{i,j}\leq 1$
and $|S_{i,j+1}-S_{i,j}|\leq 1$.
Hence, since \eqref{as:k_compact} can also be written
$$\Big|\sum\limits_{p=0}^{i-1} S_{p, j+1} - S_{p, j} \Big| \leq K,$$
it allows us to 
obtain 
\[
|\tilde G^\ep(x,y)-G^\ep(x,y)|
\leq
\ep(K +1+1/2+1/2)=\ep(K+2).
\]
We also deduce that 
\[
\big|\partial_y\tilde G^\ep(x,y)\big|=\dis\frac{\big |G^\ep(x,j\ep)-G^\ep(x,(j+1)\ep) \big|}{\ep}\leq K+2\]
while 
$$|\partial_x \tilde G^\ep(x,y)|\leq 1.$$
Moreover, we have
\[
|\tilde G^\ep(x,y)|\leq 2\ep(i  + 2) 
\]
which is bounded uniformly with respect to $\ep$ 
and $0\leq x,y\leq R<\infty$.
As a consequence of the Arzela-Ascoli theorem we deduce that, for a subsequence, $G^{\ep_n}$ converges uniformly to a continuous function $G(x,y)$ on $[0,R]\times[0,R]$ for any $0<R<\infty$.
It follows that
\[\phi^{\ep_n}(y)
\xrightarrow[\ep_n\rightarrow\infty]{}\phi(y)=
\varphi(y)-G(y,y)\varphi'(y)+\int_0^y G(x,y)\varphi''(x)\ud x
\]
uniformly on $[0,R]$.
We conclude by applying Lemma \ref{lem:G} to the function $x\mapsto G(x,y)$, with $y\geq 0$ fixed.
\qed
\begin{lemma}\label{lem:th3}
 Suppose that the discrete coefficients satisfy \fer{as:kij:strong}. Then $F^\ep$ and $G^\ep$ are uniformly bounded and converge (up to a subsequence) uniformly on compact sets.
\end{lemma}

\proof
Assumption \fer{as:kij:strong} rewrites
$$\Big| S_{i,j+1} - S_{i,j}\Big| \leq \frac{K}{j}, \qquad k_{i,j} \leq \frac{K}{j}, $$ 
so, with the same notation for $\tilde F^\ep$ as for $\tilde G^\ep$ and $\tilde k^\ep,$ we get $$|F^\ep(x,y)-\tilde F^\ep(x,y)|=\f{y-j\ep}{\ep}|F^\ep(x,(j+1)\ep)-F^\ep(x,j\ep)|\leq \f{2K}{j}\leq \f{3K}{y}\ep, $$
together (the computations are very similar to those above) with 
$$|\partial_x\tilde F^\ep|\leq \f{k_{i,j}+k_{i,{j+1}} }{\ep}\leq \f{K}{j\ep}\leq \f{K}{n_0\ep}\leq\f{2K}{x_0},$$
and 
$$|\partial_y\tilde F^\ep|\leq \f{1}{\ep}|F^\ep(x,(j+1)\ep)-F^\ep(x,j\ep)|\leq \f{3K}{y},$$
which leads to Ascoli assumptions and therefore the suitable compactness.
\qed

With such assumptions, we can take into account any $k$ of the form 
$k(x,y)dx=\frac{1}{y}k_0(x/y)dx$, including Dirac mass. If we consider such a distribution on $[0,1]$ (taken symmetric), then we can define $k_{i,j}$ as
$$k_{i,j}=k_0\Big(\Big]\f{i-1}{j-1},\f{i}{j-1}\Big[\Big)+\f{1}{2}k_0\Big(\Big\{ \f{i-1}{j-1}\Big\}\Big)+\f{1}{2}k_0\Big(\Big\{ \f{i}{j-1}\Big\}\Big)+\f{1}{2}k_0\Big(\Big\{ 0\Big\}\Big)\delta_i^1+\f{1}{2}k_0\Big(\Big\{ 0\Big\}\Big)\delta_i^{j-1}$$
with $\delta_i^j$ the Kronecker symbol.
With these notations, we have for $p\geq j-2$,
$$S_{p,j}=\sum_{i=0}^pk_{i,j}=k_0\Big(\Big[0,\f{p}{j-1}\Big[\Big)
+\f{1}{2}k_0\Big(\Big\{ \f{p}{j-1}\Big\}\Big)$$
and $S_{j-1,j}=S_{j,j}=1$, which leads to 
$$S_{p,j+1}-S_{p,j}=k_0\Big(\Big[\f{p}{j},\f{p}{j-1}\Big[\Big)+\f{1}{2}
k_0\Big(\Big\{ \f{p}{j}\Big\}\Big)-\f{1}{2}k_0
\Big(\Big\{ \f{p}{j-1}\Big\}\Big),\quad \text{if }p<j-1,$$
$$
S_{j-1,j+1}-S_{j-1,j}=-k_0\Big(\Big]\f{j-1}{j},1\Big]\Big)-\f{1}{2}
k_0\Big(\Big\{ \f{j-1}{j}\Big\}\Big),\qquad S_{j,j+1}-S_{j,j}=0,
$$
as $0\leq i\leq j$, we have for any $p\leq i$,
$$\f{p-1}{j-1}\leq \frac{p}{j},$$
the intervals $\big [\f{p}{j},\f{p}{j-1}\big[$ and $\big[\f{p-1}{j},\f{p-1}{j-1}
\big[$ are disjoint. This leads to 
$$\Big|\sum_{p=0}^i S_{p,j+1}-S_{p,j}\Big|\leq k_0\Big(\bigcup_{p=0}^i 
\Big[\f{p}{j},\f{p}{j-1}\Big[\Big)+\f{1}{2}k_0\Big(\bigcup_{p=0}^i  \Big\{\f{p}{j}\Big\}\Big)
+\f{1}{2}k_0\Big(\bigcup_{p=0}^i  \Big\{\f{p}{j-1}\Big\}\Big)\leq 2,$$
which gives the criterion \fer{as:k_compact}. The limit is then obviously given by $k(x,y)dx=\frac{1}{y}k_0(x/y)dx$.
\subsection{Discrete system}\label{app:admis}
We discuss here briefly the existence theorem for the discrete system. It is mainly an adaptation of theorem 5.1 in \cite{Laurencot-multiple}. We define the truncated system. 
Let $N>n_0$, consider the system
\begin{equation}\left\lbrace\begin{array}{l}\label{eq:discretN}
\displaystyle \frac{\ud v}{\ud t}=\lambda-\gamma v -v\dis\sum_{i=n_0}^{N-1} \tau_i u_i +2\dis\sum_{j=n_0}^N\dis\sum_{i<n_0}i k_{i,j}\beta_j u_j,\\[0.3cm]
\displaystyle \frac{\ud u_{n_0}}{\ud t}=-\mu_{n_0} u_{n_0}-\beta_{n_0}u_{n_0}-v\tau_{n_0}+2\dis\sum_{j=i+1}^N \beta_jk_{n_0,j}u_j,\qquad \text{for } ,\\[0.3cm]
\displaystyle \frac{\ud u_i}{\ud t}=-\mu_i u_i-\beta_iu_i-v(\tau_iu_i-\tau_{i-1}u_{i-1})+2\dis\sum_{j=i+1}^N \beta_jk_{i,j}u_j,\qquad \text{for } n_0<i<N ,\\[0.3cm]
\displaystyle \frac{\ud u_N}{\ud t}=-\mu_N u_N-\beta_Nu_N+v\tau_{N-1}u_{N-1}.
\end{array}\right.
\end{equation}
Existence, uniqueness and nonnegativity are immediate, we have immediately the weak formulation 
\begin{equation}\begin{array}{lll}\label{weak:4}
\dfrac{\ud }{\ud t}\bigg(v(t)\psi +\dis\sum_{i=n_0}^N u_i\varphi_i\bigg)&=&\lambda\psi-\gamma v\psi -v\dis\sum_{i=n_0}^N\mu_i u_i\varphi_i+
v\dis\sum_{i=n_0}^{N-1}\tau_iu_i(\varphi_{i+1}-\varphi_i-\psi)
\\
&&+2\dis\sum_{j=n_0+1}^N\sum_{i=n_0}^{j-1}ik_{i,j}\beta_ju_j\bigg(\frac{\varphi_i}{i}-\frac{\varphi_j}{j}\bigg)
\\
&&
+2\dis\sum_{j=n_0}^N \dis\sum_{i=1}^{n_0-1}ik_{i,j}\beta_j u_j\big(\psi-\frac{\varphi_j}{j}\big).
\end{array}
\end{equation}
Let us denote $U^N$ the infinite sequence of functions defined by $U_i^N=u_i^n$ if $n_0\leq i\leq N$, $U_i^N=0$ otherwise. The weak formulation gives moment estimates (and the moment estimates done in section 4 can then  be thought as uniform bounds on truncated systems). This model has the property of propagating moments.  

With this type of initial condition, the proof of existence is based on the Ascoli theorem for the continuous functions $U_i^N$. Thanks all the moments controlled on the initial data and the nice property of propagation of moments, we have bounds on the derivative of $v^N,U_i^N$ and therefore, we can extract convergent subsequence. The limit satisfies the equation in an integral form (see \cite{Ball-Carr-DC} for a definition). For proving uniqueness, the  procedure exposed in \cite{Ball-Carr-DC,Laurencot-multiple} applies, with a small modification due to death rates (the condition on the moment of order $1+m$ for the initial data insures the convergence of $\sum i\mu_iu_i$). 
\\

\noindent
{\bf Acknowledgments.}
We thank Pavel Dubovski for kind comments and suggestions about 
this work, and Natacha Lenuzza for fruitful discussions and help about orders of magnitude of the parameters. Part of this work has been supported by the ANR grant TOPPAZ. 
%%%%%%%%%%%%%%%%%%%%%%%%%%%%%%%%%%%%%%%%%%%%%%%%%%%%%%%%%%%%%%%%%%%%%%%%%BIBLIOGRAPHIE
%
%

%%%%%%%%%%%%%%%%%%%%%%%%%%%%%%%%%%%%%%%%%%%%%%%%%%%%%%%%%%%%%%%%%%%%%%%%%BIBLIOGRAPHIE
\bibliography{prion}
\bibliographystyle{alpha}
\end{document}